\documentclass[a4paper]{article}
\usepackage[utf8]{inputenc}
\usepackage{a4wide}
\usepackage{amsthm}
\usepackage{amsmath}
\usepackage{amsfonts}
\usepackage{amssymb,comment,xspace}
\usepackage[ruled,vlined,linesnumbered]{algorithm2e}
\usepackage{graphicx}
\usepackage{pgf}
\usepackage{tikz}
\usetikzlibrary{arrows,automata,positioning,patterns}

\usepackage{multirow}
\usepackage{mathtools}

\usepackage{url}
\usepackage[capitalize]{cleveref}

\usepackage{caption}
\usepackage{subcaption}

\newcommand{\EE}{\mathbb{E}}
\newcommand{\VV}{\mathbb{V}}
\newcommand{\PP}{\mathbb{P}}
\newcommand{\toinp}{\xrightarrow{p}}
\newcommand{\T}{\mathcal{T}}
\renewcommand{\L}{\mathcal{L}}

\newcommand{\MG}{\mathcal{M}^{>}}
\newcommand{\MS}{\mathcal{M}^{<}}

\DeclareMathOperator{\sym}{\mathfrak{S}}
\DeclareMathOperator{\weight}{w}
\DeclareMathOperator{\cyc}{cyc}
\DeclareMathOperator{\rec}{rec}
\DeclareMathOperator{\norm}{norm} % Attention, sa signification a changé par rapport au papier AofA !

\DeclareMathOperator{\desc}{desc}
\DeclareMathOperator{\inv}{inv}

\DeclareMathOperator{\bij}{\varphi}

\DeclareMathOperator{\lmax}{lmax}
\newcommand{\rfact}[2]{#1^{(#2)}}
\newcommand{\recbiased}{record-biased\xspace}

\newtheorem{theorem}{Theorem}
\newtheorem{lemma}[theorem]{Lemma}
\newtheorem{corollary}[theorem]{Corollary}
\newtheorem{proposition}[theorem]{Proposition}

\theoremstyle{definition}
\newtheorem{remark}[theorem]{Remark}

\newcommand{\ulorraine}{Université de Lorraine, CNRS, Inria, LORIA, F-54000 Nancy, France}
\newcommand{\ueiffel}{Université Gustave Eiffel, CNRS, LIGM, F-77454 Marne-la-Vallée, France}

\title{Record-biased permutations and their permuton limit}
\author{Mathilde Bouvel\thanks{\ulorraine}, Cyril Nicaud\thanks{\ueiffel}, Carine Pivoteau\thanks{\ueiffel}}

\begin{document}

\maketitle

\begin{abstract}
In this article, we study a non-uniform distribution on permutations biased by their number of records that we call \emph{record-biased permutations}. We give several generative processes for record-biased permutations, explaining also how they can be used to devise  efficient (linear) random samplers. For several classical permutation statistics, we obtain their expectation using the above generative processes, as well as their %asymptotic 
limit distributions in the regime that has a logarithmic number of records (as in the uniform case). Finally, increasing the bias to obtain a regime with an expected linear number of records, we establish the convergence of record-biased permutations to a deterministic permuton, which we fully characterize.

This model was introduced in our earlier work~\cite{AuBoNiPi2016}, in the context of realistic analysis of algorithms. We conduct here a more thorough study but with a theoretical perspective.
\end{abstract}

\section{Introduction}

Average case analysis of algorithms on discrete structures is an important research topic, at the crossroad of computer science and mathematics. It consists in stepping outside the usual framework of the worst case analysis of an algorithm by setting a distribution on the set of possible inputs of size $n$ and computing or estimating the expected running time. See~\cite{FlSe09} for a detailed account on this topic.
This approach is particularly relevant when a discrepancy from the worst case results is observed in practice: the most famous example is \textsc{QuickSort} sorting algorithm~\cite{hoare1962quicksort} whose worst case and average case running times are $\mathcal{O}(n^2)$ and $\mathcal{O}(n\log n)$, respectively, and which is implemented as the basic sorting procedure in many programming languages.

Most often in average case analysis, the uniform distribution on the set of inputs of a given size~$n$ is chosen, then letting $n$ tend to infinite. While this assumption is the most natural one from a mathematical perspective, it does not always provide an accurate view of the typical behavior of some real-life algorithms used on real-life data. For instance, some popular programming languages switched from \textsc{QuickSort} to \textsc{TimSort} or \textsc{PowerSort}~\cite{MunroW18} when they observed that the input of sorting algorithms is often already partially sorted.
To reckon with this observation and provide a more realistic analysis of the algorithms that are actually used, it is quite natural to work under the assumption of a non-uniform distribution of the inputs, chosen to account for some specific features of the data set on which the algorithm is used, while being mathematically tractable. 

The starting point of this work is the analysis, realistic in the above sense, of algorithms working on arrays of numbers, and typically sorting algorithms. In this context, the inputs can be modelled as permutations, as only the comparisons between two entries matter, not their actual values. \emph{A priori} unrelated to the computer science motivation, the study of non-uniform distributions on permutations has been an active research topic in discrete probability for several decades. Among the most studied models, one can cite the Ewens distribution~\cite{ewens1972} and the Mallows distribution~\cite{mallows1957}, both introducing an exponential bias according to a statistic on permutations (specifically, the number of cycles for the Ewens model and the number of inversions for the Mallows model). 

In our earlier work~\cite{AuBoNiPi2016}, we defined another non-uniform distribution on permutations, which we call \emph{record-biased} of parameter $\theta >0$. It also has an exponential bias: this time, the probability of a permutation~$\sigma$ is proportional to $\theta^{\rec(\sigma)}$ where $\rec(\sigma)$ is the number of records (a.k.a. left-to-right maxima) of $\sigma$. We demonstrate in~\cite{AuBoNiPi2016} that this model is meaningful in the context of presortedness, which naturally arises when studying algorithms that are designed to be efficient for almost sorted sequences~\cite{Mannila1985}. In addition, as we will review later, the record-biased distribution is the image of the Ewens distribution by the Foata bijection (or fundamental bijection~\cite[\S 10.2]{lothaire1997}). This fact can be used to prove some results on record-biased permutations: this applies to some of the results of~\cite{AuBoNiPi2016} and the present paper (although only a small proportion in both papers). 

The purpose of~\cite{AuBoNiPi2016} was to provide a fine-grained analysis of several algorithms under the record-biased distribution, in an attempt to reconcile the theoretical complexity of these algorithms with their behavior observed in practice. For example, we compared two algorithms (one naive and one more clever) solving the problem of searching both the minimum and the maximum in an array, by expressing the number of mispredictions that a branch predictor makes in each of them; this participates in explaining why the naive algorithm outerperforms in practice the clever one, while making more comparisons between entries. 
The model of record-biased permutations has also been studied by B. Corsini in~\cite{Corsini}, which describes the height of binary search trees built from record-biased permutations. 

Since then, it appeared to us that the model of record-biased permutations was likely to have nice properties, beyond those described in~\cite{AuBoNiPi2016,Corsini}. 
The purpose of the present paper is to provide a mathematical analysis of record-biased permutations, exhibiting such properties to complement these earlier results. 
An important part of our work focuses of a regime with a very strong bias towards sortedness (with $\theta$ that is linear w.r.t. the size, hence a number of records that is also linear w.r.t. the size); the interest for this regime is not only because of the nice mathematical results we can prove, but also because it is more realistic in terms of real-life applications, since using a constant $\theta$ does not change sufficiently the typical shape of the input.

A new important aspect of the present study is the description of the typical limit shape of \recbiased permutations. While our first results  describe the typical behavior of numerical quantities on \recbiased permutations, we go a step further by describing also the limit of the diagrams of \recbiased permutations, by means of their \emph{permuton} limit. Permutons have been introduced in \cite{FinitelyForcible,Permutons,PresuttiStromquist} as a permutation analogue of graphons, which allow to describe limits of dense graphs. The description of permuton limits of non-uniform random permutations has attracted quite a lot of attention lately (see \emph{e.g.} \cite[and references therein]{runsort,PHC2,PHC,Jacopo2,Jacopo3,Winkler}), and is still a developing topic. Our work also contributes to this growing  literature. 

\medskip

The rest of the article is organized as follows. 
\begin{itemize}
    \item In Section~\ref{sec:background}, we review some basics on permutations and some known results on the Ewens distribution~\cite{ewens1972}. We define record-biased permutations and present the link between Ewens permutations and record-biased permutations through Foata's bijection~\cite[\S10.2]{lothaire1997}. 
    \item In Section~\ref{sec:decision_trees}, we present several generative processes corresponding to the record-biased distribution of permutations of any given size, and how they can be used to achieve efficient random generation under this distribution. 
    \item Section~\ref{sec:stat} studies the behavior of some classical statistics on record-biased permutations. We first obtain explicit formulas for the expectations of these statistics, and study how these behave in various regimes for the parameter $\theta$, in particular when $\theta=\lambda n$ is linear w.r.t. the size $n$. 
    We next focus on the regime where $\theta$ is constant, deriving the asymptotic distributions of the studied statistics as the size goes to infinity: three are Gaussian, and the last one follows a beta distribution, asymptotically. 
    \item Finally, we describe the typical shape of the diagram of a record-biased permutation by establishing the permuton limit of random record-biased permutations of size $n$ as $n$ tends to infinity, in the linear regime where $\theta=\lambda n$. 
\end{itemize}

It should be noted that the results of Section~\ref{subsec:stat_comme_avant} already appear in~\cite{AuBoNiPi2016}, but with a different proof. While the proofs in~\cite{AuBoNiPi2016} relied on careful case analysis, the proofs in the present paper follow easily from the generative processes of Section~\ref{sec:decision_trees}. We believe that this approach is more enlightening than the former one.

%%%%%%%%%%%%%%%%%%%%%%%%%%%%%%%%%%%%%%%%%%%%%%%%%%%%%%%%%%%%%%%%%%%%%%%%%%%%%%%%
\section{Background and definition of \recbiased permutations} \label{sec:background}

%%%%%%%%%%%%%%%%%%%%%%%%%%%%%%%%%%%%%%%%
\subsection{Permutations}

For any integers~$a$ and~$b$, let $[a,b]=\{a,\ldots,b\}$ and for every integer $n\geq 1$, let $[n]=[1,n]$. 
By convention $[0]=\emptyset$. 
If $E$ is a finite set, let $\sym(E)$ denote the set of all permutations on $E$, \emph{i.e.}, of bijective maps from $E$ to itself. 
For convenience, $\sym([n])$ is written $\sym_{n}$ in the sequel. 
For a permutation~$\sigma$ in $\sym_{n}$, $n$ is called the size of~$\sigma$ and is denoted~$|\sigma|$.
 Permutations of $\sym_{n}$ can be seen in several ways (reviewed for instance in~\cite{Bona}). 
Here, we use alternately their representations as words, as diagrams, and as sets of cycles. 

A permutation $\sigma$ of $\sym_{n}$ can be represented as a word $w_1 w_2 \cdots w_n$ containing exactly once each symbol in $[n]$ 
by simply setting $w_i = \sigma(i)$ for all $i \in [n]$. 

\begin{figure}[ht]
\begin{center}
 \begin{tikzpicture}[scale=0.33]
\draw (0,0) grid (7,7);
\draw (0.5,5.5) [fill] circle (.2);
\draw (1.5,2.5) [fill] circle (.2);
\draw (2.5,1.5) [fill] circle (.2);
\draw (3.5,0.5) [fill] circle (.2);
\draw (4.5,6.5) [fill] circle (.2);
\draw (5.5,3.5) [fill] circle (.2);
\draw (6.5,4.5) [fill] circle (.2);
\end{tikzpicture}
\qquad \qquad \begin{tikzpicture}
\node (p3) at (0.5,5.5) [inner sep=1] {2};
\node (p2) at (1.5,5.5) [inner sep=1] {3};
\draw [->] (p3)  to [bend left=60] (p2);
\draw [->] (p2)  to [bend left=60] (p3);

\node (p6) at (2.5,4.25) [inner sep=1] {6};
\node (p4) at (3.0,5.25) [inner sep=1] {4};
\node (p1) at (3.5,4.25) [inner sep=1] {1};
\draw [->] (p6)  to [bend left=45] (p4);
\draw [->] (p4)  to [bend left=45] (p1);
\draw [->] (p1)  to [bend left=45] (p6);

\node (p7) at (0.5,4) [inner sep=1] {7};
\node (p5) at (1.5,4) [inner sep=1] {5};
\draw [->] (p7)  to [bend left=60] (p5);
\draw [->] (p5)  to [bend left=60] (p7);
\end{tikzpicture}
\end{center}
 \caption{The diagram and the set of cycles representations of $\tau =6321745$. \label{fig:representations}}
\end{figure}
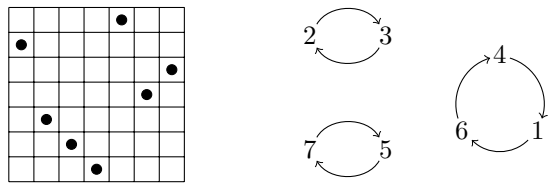

The \emph{diagram} of a permutation $\sigma$ of size $n$ consists of the set of points at coordinates $(i,\sigma(i))$ 
in an $n\times n$ grid. An example is given in~\cref{fig:representations}.

A cycle of size~$k$ in a permutation $\sigma \in \sym_{n}$ is a subset $\{i_1,\dots,i_k\}$ of $[n]$ such that 
$i_1 \stackrel{\sigma}{\mapsto} i_2 \ldots \stackrel{\sigma}{\mapsto} i_k \stackrel{\sigma}{\mapsto} i_1$. 
It is written $(i_1, i_2, \ldots, i_k)$. 
Any permutation can be decomposed as the set of its cycles. 
For instance, the cycle decomposition of $\tau$ represented by the word $6321745$ is $(32) (641) (75)$. See also \cref{fig:representations}.

%%%%%%%%%%%%%%%%%%%%%%%%%%%%%%%%%%%%%%%%
\subsection{Ewens model}

The Ewens distribution on permutations (see for instance~\cite[Ch. 4 \& 5]{Arratia} and~\cite{ewens1972} for the original article) is a generalization of the uniform distribution on $\sym_{n}$: 
the probability of a permutation depends on its number of cycles. 
Denoting $\cyc(\sigma)$ the number of cycles of any permutation $\sigma$, 
the Ewens distribution of parameter $\theta$ on $\sym_n$ (where $\theta$ is any fixed positive real number) 
gives to any $\sigma \in \sym_n$ the probability 
\begin{equation}
\frac{\theta^{\cyc(\sigma)}}{\sum_{\rho\in\sym_{n}} \theta^{\cyc(\rho)}}.
\label{eq:probEwens}
\end{equation}
(Observe that setting $\theta=1$ yields the uniform distribution on $\sym_n$.)
As seen in~\cite[Ch. 5]{Arratia}, the normalization constant $\sum_{\rho\in\sym_{n}} \theta^{\cyc(\rho)}$ is equal to the \emph{rising factorial} $\rfact\theta n$, where we recall that  for any real $x$, the rising factorial denoted $\rfact{x}{n}$ is defined by
\begin{equation}
\rfact{x}{n}=x(x + 1)\cdots(x + n-1)
\label{eq:rising}
\end{equation}
(with the convention that $\rfact{x}{0}=1$). 
Forgetting the normalizing constant, we sometimes refer to the numerator $ \theta^{\cyc(\sigma)}$ in \cref{eq:probEwens} as the \emph{weight} of the permutation $\sigma$. It is often convenient to see the weight of $\sigma$ as the product of weights associated with each of its elements. Here, this is achieved by giving weight~$\theta$ to each element that is the smallest in its cycle, and $1$ to all others.

\medskip

For fixed $\theta$, the expected number of cycles of a random permutation of size~$n$ under the Ewens distribution is $\sum_{j=0}^{n-1}\frac{\theta}{\theta+j}$ which is asymptotically equivalent to~$\theta\log n$, and it is asymptotically normal~\cite[\S4.6]{Arratia}. This will be used in the sequel, when considering the distribution on permutations studied in this article. We also rely on results on the number of weak exceedences (\emph{i.e.} indices~$i$ such that $\sigma(i)\geq i$) for the Ewens distribution, obtained in~\cite{Valentin}, as an application of general techniques for parameters defined using patterns. Namely, we shall use that the expected number of  weak exceedences is asymptotically equivalent to $n/2$ and that it is asymptotically normal.

\subsection{Generative processes for the Ewens distribution}\label{sec:generative ewens}

In the sequel, we call \emph{generative process} a way to incrementally build random permutations that follow a given distribution. They are conveniently depicted by trees, where the branches correspond to the random choices performed during the generation\footnote{This is similar to the way algorithms are translated into decision trees when establishing lower bounds in complexity, but in addition with probabilities (or weights) on the edges of the trees.}. % https://en.wikipedia.org/wiki/Decision_tree_model
As we are interested in distributions on~$\sym_n$ for every positive $n$, we consider two kinds of generative processes: either $n$ is known in advance and is used for the construction, or the generative process is in essence infinite and one halts when a size-$n$ permutation is produced.

Generative processes can directly be transformed into random generation algorithms, their complexity mainly depending on the data structures used during the process. They also prove to be very useful to establish the behavior of statistics as  the ones studied in~\cref{sec:stat}, since they describe how a random permutation can be seen as a sequence of independent random choices. We recall in this part two classical generative processes for the Ewens distribution.

\paragraph{The Chinese restaurant process.} This generative process works on permutations represented as sets of cycles, and is originally meant for uniform random permutations. A straightforward variant of the Chinese restaurant process (see~\cite[Example 2.4]{Arratia} for instance) can be used to derive any permutation according to the Ewens distribution.

One starts with an empty permutation $\sigma_0$ at step~$0$ and produces a size-$n$ random permutation~$\sigma_n$ after $n$ steps. At each step~$i\geq 1$, the new permutation $\sigma_i$ is build modifying $\sigma_{i-1}$ as follows. First, we set $\sigma_i(j)=\sigma_{i-1}(j)$ for $j\in[i-1]$. Then, we insert $i$ as described below, all the random choices being independent. 
\begin{itemize}
    \item With probability~$\frac{\theta}{\theta + i-1}$, we create a new cycle containing~$i$ only, {\em i.e.}, we set $\sigma_i(i)=i$.
    \item With probability~$\frac{i-1}{\theta + i-1}$, we insert~$i$ in one of the $i-1$ positions available in an existing cycle, these positions being equally likely. 
    Equivalently, for each $j\in[i-1]$, with probability~$\frac{1}{\theta + i-1}$, set $\sigma_i(i)=j$ and $\sigma_i(\sigma_{i-1}^{-1}(j))=i$ (which corresponds to inserting $i$ before $j$ in its cycle). 
\end{itemize}
This process is infinite and does not depend on~$n$. However, stopping after~$n$ steps produces a random permutation of size~$n$ according to the Ewens distribution of parameter~$\theta$. Indeed, with the notation $c_{\sigma}(i) = 1$ if the $i$ is the minimum in its cycle in $\sigma$ and $c_{\sigma}(i) = 0$ otherwise, it generates the permutation $\sigma \in \sym_n$ with probability 
\[
\prod_{i=1}^{n}\frac{\theta^{c_{\sigma}(i)}}{i-1+\theta} = \frac{\theta^{\cyc(\sigma)}}{\theta^{(n)}}.
\] 

This process can be represented by a tree, which is also infinite in this instance.  Its root, at height~$0$, is the empty permutation and it has only one child which is a single cycle containing~$1$. Then, each node at height~$i$ has $i+1$ children obtained as described above. 

Trimming the tree at height~$n$ creates leaves that are all the permutations of size~$n$.
For example, the tree for $\sym_{3}$ depicted on \cref{fig:chinese-tree}. 
In this figure (as well as in similar figures later), the weight of a permutation corresponding to a leaf of the tree is shown in red, to the right of the permutation. As explained earlier, this weight can be obtained as the product of weights associated with its elements. The weight associated with an element is shown in blue on the figures, on the edge of the tree that correspond to the insertion of the considered element.

\begin{figure}[t]
\centering
\begin{minipage}{.48\textwidth}
\vspace{-7mm}
\begin{tikzpicture}[style={scale=0.78}]
  \node {$\emptyset$} [level distance=10mm, minimum size=10mm, grow=right]
    child { [level distance=15mm, sibling distance=41mm, minimum height=12mm] node (1) {}
        child { [level distance=36mm, sibling distance=16mm , minimum size=17mm] 
          node (21) {} 
          child {node (312) {} edge from parent node[near end,yshift=2mm,color=blue] {$1$}} 
          child {node (231) {} edge from parent node[near end,yshift=2mm,color=blue] {$1$}}
          child {node (213) {} edge from parent node[near end,yshift=2mm,color=blue] {$\theta$}}
        }
        child { [level distance=36mm, sibling distance=17mm,inner sep=10pt, minimum size=17mm] 
          node (12) {}
          child {node (321) {} edge from parent node[near end,yshift=2mm,color=blue] {$1$}}
          child {node (132) {} edge from parent node[near end,yshift=2mm,color=blue] {$1$}}
          child {node (123) {} edge from parent node[near end,yshift=2mm,color=blue] {$\theta$}}
        }
      };     

    \node at (1) [xshift=-2mm,yshift=2mm,color=blue] {$\theta$};
    \node at (21) [xshift=-4mm,yshift=2mm,color=blue] {$1$};
    \node at (12) [xshift=-4mm,yshift=-2mm,color=blue] {$\theta$};

    \node (p1) at (1) [xshift=2mm, inner sep=2] {1};
    \draw [->] (p1) [xshift=1.5mm,yshift=1mm] arc(150:-130:0.25);

    \node (p21) at (21) [xshift=1.5mm, inner sep=2] {1};
    \node (p21-2) [inner sep=2,right=of p21,xshift=-7.5mm] {2};
    \draw [->] (p21)  to [bend left=65] (p21-2);
    \draw [->] (p21-2)  to [bend left=65] (p21);

    \node (p12) at (12) [yshift=3mm,xshift=2mm, inner sep=3] {1};
    \node (p12-2) [below=of p12,yshift=9mm, inner sep=3] {2};
    \draw [->] (p12)[xshift=1.5mm,yshift=1mm] arc(150:-130:0.25);
    \draw [->] (p12-2)[xshift=1.5mm,yshift=1mm] arc(150:-130:0.25);

    \node (p312) at (312) [xshift=-5mm,yshift=4mm, inner sep=1] {1};
    \node (p312-2) [inner sep=1,right=of p312,xshift=-7mm] {2};
    \node (p312-3) [inner sep=1,below=of p312,xshift=3mm, yshift=9mm] {3};
    \draw [->] (p312)  to [bend left=45] (p312-2);
    \draw [->] (p312-2)  to [bend left=45] (p312-3);
    \draw [->] (p312-3)  to [bend left=45] (p312);

    \node (p231) at (231) [xshift=-5mm, yshift=2mm,inner sep=1] {1};
    \node (p231-2) [inner sep=1,right=of p231,xshift=-7mm] {3};
    \node (p231-3) [inner sep=1,below=of p231,xshift=3mm, yshift=9mm] {2};
    \draw [->] (p231)  to [bend left=45] (p231-2);
    \draw [->] (p231-2)  to [bend left=45] (p231-3);
    \draw [->] (p231-3)  to [bend left=45] (p231);

    \node (p213) at (213) [yshift=-2mm,xshift=-5mm, inner sep=2] {1};
    \node (p213-2) [inner sep=2,right=of p213,xshift=-7.5mm] {2};
    \node (p213-3) [right=of p213-2, xshift=-7mm, inner sep=3] {3};
    \draw [->] (p213)  to [bend left=65] (p213-2);
    \draw [->] (p213-2)  to [bend left=65] (p213);
    \draw [->] (p213-3)[xshift=1.5mm,yshift=1mm] arc(150:-130:0.25);

    \node (p321) at (321) [yshift=4mm,xshift=-5mm, inner sep=3] {1};
    \node (p321-2) [inner sep=1,right=of p321,xshift=-7.5mm] {3};
    \node (p321-3) [right=of p321-2, xshift=-7mm, inner sep=3] {2};
    \draw [->] (p321)  to [bend left=65] (p321-2);
    \draw [->] (p321-2)  to [bend left=65] (p321);
    \draw [->] (p321-3)[xshift=1.5mm,yshift=1mm] arc(150:-130:0.25);

    \node (p132) at (132) [yshift=1mm,xshift=-5mm, inner sep=3] {3};
    \node (p132-2) [inner sep=1,right=of p132,xshift=-7.5mm] {2};
    \node (p132-3) [right=of p132-2, xshift=-7mm, inner sep=3] {1};
    \draw [->] (p132)  to [bend left=65] (p132-2);
    \draw [->] (p132-2)  to [bend left=65] (p132);
    \draw [->] (p132-3)[xshift=1.5mm,yshift=1mm] arc(150:-130:0.25);

    \node (p123) at (123) [yshift=-3mm,xshift=-5mm, inner sep=3] {1};
    \draw [->] (p123)[xshift=1.5mm,yshift=1mm] arc(150:-130:0.25);
    \node (p123-2) [right=of p123,xshift=-6mm, inner sep=3] {2};
    \draw [->] (p123-2)[xshift=1.5mm,yshift=1mm] arc(150:-130:0.25);
    \node (p123-3) [right=of p123-2,xshift=-6mm, inner sep=3] {3}; 
    \draw [->] (p123-3)[xshift=1.5mm,yshift=1mm] arc(150:-130:0.25);

    \node at (p123) [xshift=24mm,color=red] {\small $\theta^3$};
    \node at (p132) [xshift=24mm,color=red]  {\small $\theta^2$};
    \node at (p213) [xshift=24mm,color=red]  {\small $\theta^2$};
    \node at (p231) [xshift=23.5mm,yshift=-1mm,color=red]{\small $\theta$};
    \node at (p321) [xshift=24mm,color=red]  {\small $\theta^2$};
    \node at (p312) [xshift=23.5mm,yshift=-1mm,color=red]{\small $\theta$};
\end{tikzpicture}
\vspace*{-12.5mm}
\caption{the Chinese Restaurant process for permutations in $\sym_{3}$ in the Ewens model.}\label{fig:chinese-tree}
\end{minipage}
\hfill
\begin{minipage}{.47\textwidth}
\vspace{1mm}
\begin{tikzpicture}[style={scale=0.78}]
  \node {$\emptyset$} [level distance=12mm, grow=right]
    child { [level distance=13mm, sibling distance=27mm] node (1) [fill=white,anchor=west] {$(1$}
        child { [level distance=28mm, sibling distance=13mm] 
          node (12) [fill=white,anchor=west] {$(1)(2$} 
          child {node (123) [fill=white,anchor=west] {$(1)(2)(3)$} edge from parent node[below,color=blue] {$~\theta$}} 
          child {node (132) [fill=white,anchor=west] {$(1)(23)$} edge from parent node[above,color=blue] {$1$}}
          edge from parent node[above,xshift=1mm,color=blue] {$\theta$}
        }
        child { [level distance=28mm, sibling distance=13mm] 
          node (21) [fill=white,anchor=west] {$(12\phantom{)(}$}
          child {node (213) [fill=white,anchor=west] {$(12)(3)$} edge from parent node[below,color=blue] {$~\theta$}}
          child {node (231) [fill=white,anchor=west] {$(123)$} edge from parent node[above,color=blue] {$1$}}
          edge from parent node[above,color=blue] {$1$}
        }
        child { [level distance=28mm, sibling distance=14mm] 
          node (31) [fill=white,anchor=west] {$(13\phantom{)(}$}
          child {node (321) [fill=white,anchor=west] {$(13)(2)$} edge from parent node[below,color=blue] {$~\theta$}}
          child {node (312) [fill=white,anchor=west] {$(132)$} edge from parent node[above,color=blue] {$1$}}
          edge from parent node[above,xshift=0.5mm,yshift=2.5mm,color=blue] {$1$}
        }
        edge from parent node[above,color=blue] {$\theta$}
      };     

    \node at (123) [xshift=9.5mm,color=red] {\small $\theta^3$};
    \node at (132) [xshift=11mm,color=red]  {\small $\theta^2$};
    \node at (213) [xshift=11mm,color=red]  {\small $\theta^2$};
    \node at (231) [xshift=11.5mm,color=red]{\small $\theta$};
    \node at (321) [xshift=11mm,color=red]  {\small $\theta^2$};
    \node at (312) [xshift=11.5mm,color=red]{\small $\theta$};
\end{tikzpicture}
\caption{Feller coupling for permutations in $\sym_{3}$ in the Ewens model.}\label{fig:fellers-tree}
\end{minipage}
\end{figure}
\paragraph{The Feller coupling.}

The {\em Feller coupling} is another classical generative process that can be used for the Ewens distribution (see~\cite[p.16]{Arratia}, for instance). Contrary to the Chinese restaurant process, the size $n$ has to be known in advance, and there is one (finite) tree associated with each $n$. This generative process can be described as follows. At any given step, except for the initialization, the partial representation of the permutation consists of a set of cycles and a sequence, called the \emph{open cycle}. The open cycle represents the current cycle under construction, and the other cycles are complete: they remain untouched in the sequel and will be cycles of the generated permutation at the end. At each step~$i$,  we perform an (independent) random choice as follows:
\begin{itemize}
    \item with probability $\frac{\theta}{\theta+n-i}$, the open cycle is closed, adding the newly formed cycle to the list, and we start a new open cycle containing the smallest unused value;
    \item for every unused value $j\in[n]$, with probability $\frac{1}{\theta+n-i}$ we add $j$ at the end of the open cycle.
\end{itemize}
We start with no cycle and no open cycle, so the first step deterministically produces a open cycle containing  the value $1$ only.

For the tree representation depicted in \cref{fig:fellers-tree}, we write the cycles in their order of generation, which naturally orders them by increasing values of their smallest elements. The open cycle is written with only an open parenthesis: $(13$ represents the open cycle $1\rightarrow3\rightarrow$. On each edge, we put a weight $\theta$ or $1$ depending on the chosen case, to match with the probabilities after normalization. One can readily check that the probability that a given size-$n$ permutation $\sigma$ generated by this process is:

\[
\prod_{i=1}^{n}\frac{\theta^{c_{\sigma}(i)}}{n-i+\theta} = \frac{\theta^{\cyc(\sigma)}}{\theta^{(n)}}.
\] 
with $c_{\sigma}(i) = 1$ if  $i$ is the minimum in its cycle in $\sigma$ and $c_{\sigma}(i) = 0$ otherwise.

%%%%%%%%%%%%%%%%%%%%%%%%%%%%%%%%%%%%%%%%
\subsection{Foata's bijection} \label{sec:Foata}

For a permutation $\sigma \in \sym_{n}$, and for $i\in[n]$, we say that 
there is a \emph{record} at position~$i$ in~$\sigma$ (and subsequently, that $\sigma(i)$ is a record) if $\sigma(i)>\sigma(j)$ for every $j\in [i-1]$. 
Records are also called \emph{max-records} or \emph{left-to-right maxima} in other works. 

\emph{Foata's bijection}~\cite[\S10.2]{lothaire1997}, also called the \emph{fundamental bijection}, is a bijection from $\sym_n$ to $\sym_n$ that maps the number of cycles to the number of records in a permutation. 
This bijection $\bij$ is the following transformation. 
\begin{enumerate}
\setlength\itemsep{-1mm}
\item Given $\sigma$ a permutation of size $n$, consider the cycle decomposition of $\sigma$. 
\item Write every cycle starting with its maximal element, 
and write the cycles in increasing order of their maximal (\emph{i.e.,} first) element. 
\item Erase the parenthesis to get the word representation of $\bij(\sigma)$. 
\end{enumerate}
For instance $\bij\big( (57)(32)(416) \big) = 3264175$. 

This transformation is a bijection, and clearly satisfies that each cycle of $\sigma$ corresponds to one record of $\bij(\sigma)$, and conversely. 
For references and details about this bijection, see for example~\cite[p. 109--110]{Bona} or~\cite[\S10.2]{lothaire1997}.

\subsection{Our model: Record-biased permutations}

Recall the definition of a record, which is essential here:
for $i\in [n]$, $\sigma(i)$ is a record of~$\sigma \in \sym_{n}$ if $\sigma(i)>\sigma(j)$ for every $j\in [i-1]$. 
In the word representation of permutations, records are therefore elements that have no larger elements to their left. 
And in the diagrams of permutations, records are elements that have no higher elements to their left.
An element that is not a record is called a \emph{non-record}.
We denote by $\rec(\sigma)$ the number of records of a permutation $\sigma$.

In a similar fashion as for Ewens distribution, we consider the non-uniform distribution of parameter~$\theta$ that gives the probability
\[
\frac{\theta^{\rec(\sigma)}}{\sum_{\rho\in\sym_{n}} \theta^{\rec(\rho)}}
\] 
to any permutation~$\sigma \in \sym_{n}$. 
We call them {\em \recbiased} permutations. 

We define the weight of a permutation $\sigma$ as $\weight(\sigma)=\theta^{\rec(\sigma)}$. 
Like for the Ewens distribution, the normalization constant $\sum_{\rho\in\sym_{n}} \theta^{\rec(\rho)}$ is equal to the rising factorial $\rfact\theta n$ (it is for instance a direct consequence of Foata's bijection).

Foata's bijection provides a very tight connection between Ewens permutations and \recbiased permutations. Indeed, $\sigma$ is a random Ewens permutation in $\sym_n$ if and only $\bij(\sigma)$ is a random \recbiased permutation in $\sym_n$. This simple fact will be useful on several occasions in \cref{sec:decision_trees,sec:stat}.

%%%%%%%%%%%%%%%%%%%%%%%%%%%%%%%%%%%%%%%%%%%%%%%%%%%%%%%%%%%%%%%%%%%%%%%%%%%%%%%%
\section{Generative processes for \recbiased permutations} \label{sec:decision_trees}

Similarly to Section~\ref{sec:generative ewens}, we propose several generative processes for \recbiased permutations. They will be used in the sequel to compute the expected values and some limit law for several statistics of interest.

\newcommand{\tofill}[3]{
\setlength{\tabcolsep}{3pt}
\begin{tabular}{|p{4pt}|p{4pt}|p{4pt}|}
\hline#1&#2&#3\\\hline
\end{tabular}}

\begin{figure}[t]
\centering
\begin{minipage}{.4\textwidth}
\vspace*{6.5mm}
\begin{tikzpicture}[style={scale=0.78}]
  \node {$\emptyset$} [level distance=10mm, grow=right]
    child { [level distance=12mm, sibling distance=28.5mm] node (1) [fill=white,anchor=west] {$[3$}
        child { [level distance=17mm, sibling distance=13mm] 
          node (12) [fill=white,anchor=west] {$[2~[3]$} 
          child {node (123) [fill=white,anchor=west] {$[1]~[2]~[3]$} edge from parent node[below,color=blue] {$~\theta$}} 
          child {node (132) [fill=white,anchor=west] {$[21]~[3]$} edge from parent node[above,color=blue] {$1$}}
          edge from parent node[above,xshift=1mm,color=blue] {$\theta$}
        }
        child { [level distance=17mm, sibling distance=13mm] 
          node (21) [fill=white,anchor=west] {$~[32$}
          child {node (213) [fill=white,anchor=west] {$[1]~[32]$} edge from parent node[below,color=blue] {$~\theta$}}
          child {node (231) [fill=white,anchor=west] {$[321]$} edge from parent node[above,color=blue] {$1$}}
          edge from parent node[above,color=blue] {$1$}
        }
        child { [level distance=17mm, sibling distance=13mm] 
          node (31) [fill=white,anchor=west] {$~[31$}
          child {node (321) [fill=white,anchor=west] {$[2]~[31]$} edge from parent node[below,color=blue] {$~\theta$}}
          child {node (312) [fill=white,anchor=west] {$[312]$} edge from parent node[above,color=blue] {$1$}}
          edge from parent node[above,xshift=0.5mm,yshift=2.5mm,color=blue] {$1$}
        }
        edge from parent node[above,color=blue] {$\theta$}
      };     

    \node at (123) [xshift=10mm,color=red] {\small $\theta^3$};
    \node at (132) [xshift=11mm,color=red]  {\small $\theta^2$};
    \node at (213) [xshift=12mm,color=red]  {\small $\theta^2$};
    \node at (231) [xshift=12mm,color=red]{\small $\theta$};
    \node at (321) [xshift=12mm,color=red]  {\small $\theta^2$};
    \node at (312) [xshift=12mm,color=red]{\small $\theta$};
\end{tikzpicture}
\caption{Generative process for \recbiased permutations in $\sym_{3}$, viewed as sequences of sequences.}\label{fig:subsequences-tree}
\end{minipage}
\hfill
\begin{minipage}{.56\textwidth}
\begin{tikzpicture}[style={scale=0.78}]
  \node (xxx){\tofill{}{}{}} [level distance=26mm,sibling distance=28mm, grow=right]
      child { [sibling distance=15mm] 
        node (xx1) [fill=white] {\tofill{}{}{1}} 
        child { node (x21) [fill=white] {\tofill{}{2}{1}} 
                child { node (321) [fill=white] {\tofill{3}{2}{1}} 
                        edge from parent node[above,color=blue] {$\theta$}}
                edge from parent node[below, color=blue] {$1$}
        } 
        child { node (2x1) [fill=white] {\tofill{2}{}{1}} 
                child { node (231) [fill=white] {\tofill{2}{3}{1}} 
                        edge from parent node[above,color=blue] {$\theta$}}
                edge from parent node[above, color=blue] {$\theta$}
        }
        edge from parent node[below, xshift=-1mm, color=blue] {$1$}
      }
      child { [sibling distance=15mm] 
        node (x1x) [fill=white] {\tofill{}{1}{}} 
        child { node (x12) [fill=white] {\tofill{}{1}{2}} 
                child { node (312) [fill=white] {\tofill{3}{1}{2}} 
                        edge from parent node[above,color=blue] {$\theta$}}
                edge from parent node[below, color=blue] {$1$}
        } 
        child { node (21x) [fill=white] {\tofill{2}1{}} 
                child { node (213) [fill=white] {\tofill{2}{1}{3}} 
                        edge from parent node[above,color=blue] {$\theta$}}
                edge from parent node[above, color=blue] {$\theta$}
        }
        edge from parent node[above, xshift=-4mm, color=blue] {$1$}
      }
      child { [sibling distance=15mm] 
        node (1xx) [fill=white] {\tofill{1}{}{}} 
        child { node (1x2) [fill=white] {\tofill{1}{}{2}} 
                child { node (132) [fill=white] {\tofill{1}{3}{2}} 
                        edge from parent node[above,color=blue] {$\theta$}}
                edge from parent node[below, color=blue] {$1$}
        }
        child { node (12x) [fill=white] {\tofill{1}{2}{}} 
                child { node (123) [fill=white] {\tofill{1}{2}{3}} 
                        edge from parent node[above,color=blue] {$\theta$}}
                edge from parent node[above, color=blue] {$\theta$}
        } 
        edge from parent node[above, xshift=-1mm, color=blue] {$\theta$}
      }
    ;     

    \node at (123) [xshift=10mm,color=red] {\small $\theta^3$};
    \node at (132) [xshift=10mm,color=red] {\small $\theta^2$};
    \node at (213) [xshift=10mm,color=red] {\small $\theta^2$};
    \node at (231) [xshift=10mm,color=red] {\small $\theta^2$};
    \node at (321) [xshift=10mm,color=red] {\small $\theta$};
    \node at (312) [xshift=10mm,color=red] {\small $\theta$};
\end{tikzpicture}
\caption{Generative process for \recbiased permutations in $\sym_{3}$, viewed as words.}\label{fig:array-tree}
\end{minipage}
\end{figure}

%%%%%%%%%%%%%%%%%%%%%%%%%%%%%%%%%%%%%%%%
\subsection{Two generative processes for the word representation}\label{sec:process_words}

The {\em Feller coupling} presented above can readily be transformed into a generative process for \recbiased permutations. The modifications are the following:
\begin{itemize}
    \item we work directly on sequences of values instead of cycles, which are concatenated at the end;
    \item when a new open sequence is created, it contains the largest available value instead of the smallest one;
    \item the sequences are created from right to left.
\end{itemize}
This construction directly ensures that the largest value of each sequence is a record. The tree associated with this generative process on $\sym_3$ is depicted in Fig.~\ref{fig:subsequences-tree}.

\medskip

The tree depicted in~\cref{fig:array-tree} represents a different generative process for \recbiased permutations viewed as sequences. Given the target size $n$, the idea is to start with a sequence made of $n$ empty slots, and then to incrementally place the values ranging from $1$ to $n$, one by one. At each step~$i$ an independent random choice is performed:
\begin{itemize}
    \item with probability $\frac{\theta}{\theta+n-i}$, the value $i$ is placed at the leftmost empty slot;
    \item for every index $j$ in the array which corresponds to an empty slot but the leftmost one, $i$ is placed at the $j$-th slot with probability $\frac1{\theta+n-i}$.
\end{itemize}
The construction is readily seen to associate with each permutation $\sigma$ the probability $\theta^{\rec(\sigma)}/\theta^{(n)}$ which characterizes \recbiased permutations.

%%%%%%%%%%%%%%%%%%%%%%%%%%%%%%%%%%%%%%%%
\subsection{One generative process for diagrams}\label{sec:process_diagrams}

\definecolor{green3}{RGB}{0, 160, 150}
\definecolor{green2}{RGB}{0, 110, 100}
\definecolor{green1}{RGB}{0, 60, 50}

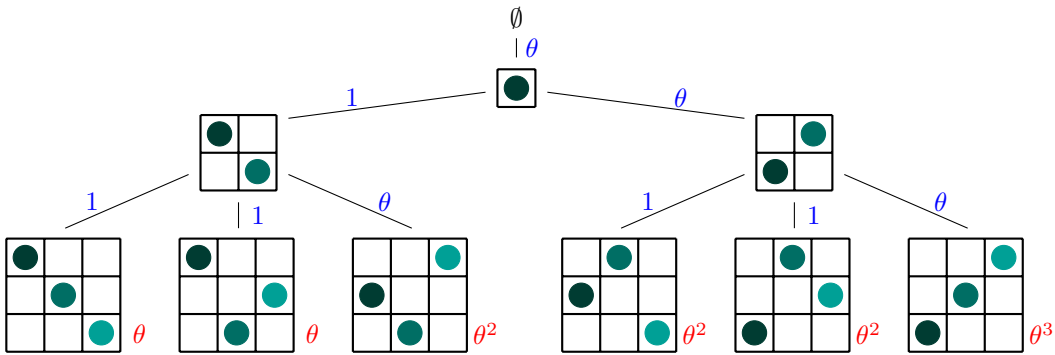
\begin{figure}[t]
%\centering
\hspace*{-4mm}
\begin{tikzpicture}[style={scale=0.4},
    complexnode/.pic={
\draw[step=.5cm,thick](-1,-1) grid (.5,.5);}]]
  \node {$\emptyset$} [level distance=13mm, minimum size=10mm]
    child { [level distance=10mm, sibling distance=150mm, minimum size=8mm,anchor=north] 
    node (1) {  
        \begin{tikzpicture}
                \draw[step=.5cm,thick](0,0) grid (.5,.5);
                \draw (0.25,.25) [fill,color=green1] circle (.16);
        \end{tikzpicture}}
        child { [level distance=25mm, sibling distance=57mm , minimum size=13mm,anchor=north east] 
        node (21) [yshift=2mm] {   
                \begin{tikzpicture}
                        \draw[step=.5cm,thick](0,0) grid (1,1);
                        \draw (0.25,.75) [fill,color=green1] circle (.16);
                        \draw (0.75,.25) [fill,color=green2] circle (.16);
                \end{tikzpicture}}
          child {[level distance=7mm,anchor=north] node (321) {   
                \begin{tikzpicture}
                        \draw[step=.5cm,thick](0,0) grid (1.5,1.5);
                        \draw (0.25,1.25) [fill,color=green1] circle (.16);
                        \draw (0.75,0.75) [fill,color=green2] circle (.16);
                        \draw (1.25,0.25) [fill,color=green3] circle (.16);
                        \node at (1.75,.9) [color=red] {$\theta$};
                        \node at (-.2,.9) [color=white] {$x$};
                \end{tikzpicture}}
                edge from parent node[left,xshift=2mm,color=blue] {$1$}} 
          child {[level distance=7mm,anchor=north] node (312) {  
                \begin{tikzpicture}
                        \draw[step=.5cm,thick](0,0) grid (1.5,1.5);
                        \draw (0.25,1.25) [fill,color=green1] circle (.16);
                        \draw (0.75,0.25) [fill,color=green2] circle (.16);
                        \draw (1.25,0.75) [fill,color=green3] circle (.16);
                        \node at (1.75,.9) [color=red] {$\theta$};
                        \node at (-.2,.9) [color=white] {$x$};
                \end{tikzpicture}}
                edge from parent node[right,xshift=-4mm,color=blue] {$1$}}
          child {[level distance=7mm,anchor=north] node (213) {  
                \begin{tikzpicture}
                        \draw[step=.5cm,thick](0,0) grid (1.5,1.5);
                        \draw (0.25,0.75) [fill,color=green1] circle (.16);
                        \draw (0.75,0.25) [fill,color=green2] circle (.16);
                        \draw (1.25,1.25) [fill,color=green3] circle (.16);
                        \node at (1.75,.9) [color=red] {$\theta^2$};
                        \node at (-.2,.9) [color=white] {$x$};
                \end{tikzpicture}} 
                edge from parent node[right,xshift=-2mm,color=blue] {$\theta$}}edge from parent node[left,xshift=2mm,yshift=1mm,color=blue] {$1$}
        }
        child { [level distance=25mm, sibling distance=57mm , minimum size=13mm,anchor=north west]  
        node (12)[yshift=2mm] {  
                \begin{tikzpicture}
                        \draw[step=.5cm,thick](0,0) grid (1,1);
                        \draw (0.25,.25) [fill,color=green1] circle (.16);
                        \draw (0.75,.75) [fill,color=green2] circle (.16);
                \end{tikzpicture}}
          child {[level distance=7mm,anchor=north] node (231) {  
                \begin{tikzpicture}
                        \draw[step=.5cm,thick](0,0) grid (1.5,1.5);
                        \draw (0.25,0.75) [fill,color=green1] circle (.16);
                        \draw (0.75,1.25) [fill,color=green2] circle (.16);
                        \draw (1.25,0.25) [fill,color=green3] circle (.16);
                        \node at (1.75,.9) [color=red] {$\theta^2$};
                        \node at (-.2,.9) [color=white] {$x$};
                \end{tikzpicture}} 
                edge from parent node[left,xshift=2mm,color=blue] {$1$}}
          child {[level distance=7mm,anchor=north] node (132) {  
                \begin{tikzpicture}
                        \draw[step=.5cm,thick](0,0) grid (1.5,1.5);
                        \draw (0.25,0.25) [fill,color=green1] circle (.16);
                        \draw (0.75,1.25) [fill,color=green2] circle (.16);
                        \draw (1.25,0.75) [fill,color=green3] circle (.16);
                        \node at (1.75,.9) [color=red] {$\theta^2$};
                        \node at (-.2,.9) [color=white] {$x$};
                \end{tikzpicture}} 
                edge from parent node[right,xshift=-4mm,color=blue] {$1$}}
          child {[level distance=7mm,anchor=north] node (123) {  
                \begin{tikzpicture}
                        \draw[step=.5cm,thick](0,0) grid (1.5,1.5);
                        \draw (0.25,0.25) [fill,color=green1] circle (.16);
                        \draw (0.75,0.75) [fill,color=green2] circle (.16);
                        \draw (1.25,1.25) [fill,color=green3] circle (.16);
                        \node at (1.75,.9) [color=red] {$\theta^3$};
                        \node at (-.2,.9) [color=white] {$x$};
                \end{tikzpicture}}
                edge from parent node[right,xshift=-2mm,color=blue] {$\theta$}}
          edge from parent node[right,xshift=-2mm,yshift=1mm,color=blue] {$\theta$}
        }
        edge from parent node[right,xshift=-2mm,color=blue] {$\theta$}
      };     
\end{tikzpicture}
\vspace*{-8mm}
\caption{Generative process for \recbiased permutations in $\sym_{3}$, viewed as diagrams.}\label{fig:diagram-tree}
\end{figure}

For permutations viewed as diagrams, one can also design a generative process, starting with an empty diagram and inserting points from left to right. Note that this generative process does not require the target size $n$ to be known in advance. This process will prove quite useful in the sequel to compute the expected values of some of our chosen statistics. 

The root of the tree is the empty diagram.
At each step~$i$, we add a new row and a new column in the diagram, and we put the new point at their crossing. 
The added column will always be the rightmost ({\em i.e.,} the $i$-th column). The placement of the point in the newly added column is chosen randomly and independently at each step~$i$:
\begin{itemize}
    \item it is on the highest row with probability $\frac{\theta}{\theta+i-1}$;
    \item every other possibility has probability
    $\frac{1}{\theta+i-1}$.
\end{itemize}
The associated tree is depicted on Fig.~\ref{fig:diagram-tree}, where the weights on the edges are such that the total weight of a permutation $\sigma$, \emph{i.e.} the product of the weights along the path from the root to the diagram of $\sigma$, is precisely $\sigma^{\rec(\sigma)}$.

%%%%%%%%%%%%%%%%%%%%%%%%%%%%%%%%%%%%%%%%

\subsection{Random generation}\label{sec:random generation}

As mentioned in~\cite[\S 2.1]{Valentin}, one can easily obtain a linear time and space algorithm to generate a random permutation according to the classical Ewens distribution, 
using the variant of the Chinese restaurant process reviewed above. 
Using the {\em fundamental bijection} of Foata, this gives us a linear random sampler for permutations according to the \recbiased distribution. 
This was fully reviewed in our previous work~\cite{AuBoNiPi2016}. 

Building on the same ideas, we can design linear samplers that generate directly \recbiased permutations based on our three generative processes, provided there exists an algorithm that allows to perform each step in constant time. We assume in the sequel that we can sample a uniform real number in $[0,1]$ in constant time, which implies that we can also generate a random integer of $\{1,\ldots,n\}$ in constant time.\footnote{An analysis for a more accurate complexity model is doable, \emph{e.g.} considering that we can only generate uniform random bits in constant time, but it is not the kind of study we want to achieve in the present article.}

\subsubsection*{Sequence representation}

To get a linear time sampler following the process of \cref{fig:array-tree}, we maintain a set $S$, initially the interval $[n]$, which represents the set of positions still empty in the array $\sigma$; 
and we need to be able to perform the following actions in constant time: 
\begin{itemize}
 \item[(i)] remove the minimum from $S$, which  will be needed when inserting a record in $\sigma$; 
 \item[(ii)] choose uniformly at random one value different from the minimum in $S$, and remove it from~$S$, which will be needed  when inserting a non-record in $\sigma$.  
\end{itemize}
This is achieved by maintaining three linear-space data-structures together: a linked list~$L$, and  two arrays~$A$ and $invA$. 
Both $L$ and $A$ store the positions of $\sigma$ which still need to be filled, \emph{i.e.} the set $S$ (but with different data-structures, which allow for different operations to be performed efficiently). 
Then, $invA$ represents the functional inverse of $A$ and allows to find in constant time the position of any given value in $A$. 

The data-structure $L$ is a doubly linked list of the values in the set $S$, in increasing order. Its cells are kept in an array such that for any value~$i$, the cell containing~$i$ is found at index~$i$, which gives direct access to each value. 
This allows to easily find the minimum, which is always the first value in the list. 
And each time a value is removed from the set $S$ (being the minimum or not), it can be removed from the list $L$ by updating pointers in constant time. 

The array~$A$ records all the positions still empty in $\sigma$, except the minimal one. 
Initially, $A$ is of length $n-1$, with $A[i]=i+1$. 
After elements $1, 2 \dots, i$ have been inserted in $\sigma$, $A$ is of length $n-i-1$ and contains the values of the set $S$ (\emph{i.e.} the still empty positions in~$\sigma$) except the minimal one. 
Having~$S$ as an array without gaps allows us to choose the position of a non-record uniformly at random in constant time. 
To maintain this array, each time we insert an element as in (ii), we choose uniformly at random an index~$k$ in the array~$A$ (the corresponding position~$A[k]$ in $\sigma$ is now filled). To take this position out of~$A$, we swap $A[k]$ with the last element of $A$, and then we shorten $A$ by one. 
When inserting an element as in (i), one of the positions recorded in~$A$ becomes the new minimum of the still empty positions in $\sigma$ and must be taken out from $A$.
The value of this element is known using $L$, and we can find its index~$k$ in $A$ in constant time using~$invA$. We take it out from~$A$ as in the previous case.

\medskip

The generative process of \cref{fig:subsequences-tree} can be transformed into a linear sampler for \recbiased permutations, using almost the same data structure as above, replacing the minimum by the maximum.
The only slight difference is that the maximum of the remaining values can be picked for action~(ii). 
This is easy to implement by keeping the maximum in $A$ and removing it (instead of the next maximum) from $A$ during action (i).
The sequence of sub-sequences is a simple list of lists in which insertions can be performed in constant time. Flattening the list in the end yields another linear time sampler for \recbiased permutations. 

\subsubsection*{Diagram representation}
To get a linear time sampler following the process of \cref{sec:process_diagrams}, we first note that choosing the height of the new row for a non-record is equivalent to choosing uniformly at random a column~$j$ and placing the new row just under the point in column $j$. 
Thus, it is enough to keep at any time a linked list of the indices of the columns ordered by decreasing height of the point they contain, starting with the highest. 
When a non-record point is added, its height is determined by choosing a column in $[i-1]$. 
And, at each step, the insertion of the index~$i$ of the last column into the list can also be done in constant time, provided we use a linked list with direct access to its cells. 
This is achieved, as before, by storing the cells in an array such that the cell containing column~$i$ is found at index~$i$.
Viewed as a word, the resulting list 
(in reverse order) corresponds to~$\sigma^{-1}$, from which we can recover $\sigma$ in linear time.

%%%%%%%%%%%%%%%%%%%%%%%%%%%%%%%%%%%%%%%%%%%%%%%%%%%%%%%%%%%%%%%%%%%%%%%%%%%%%%%%
\section{Behavior of some classical statistics} \label{sec:stat}

In this section, we study the behavior of several classical statistics on \recbiased permutations with parameter $\theta$ (which may be a fixed positive real number, or a function of the size~$n$). 
Some of the results below were already present in~\cite{AuBoNiPi2016}: mostly, until the expectation of the considered statistics. 
In~\cite{AuBoNiPi2016}, further results were then established: they concern the behavior of these expectations in various regimes for $\theta$, 
and some applications to the analysis of algorithms (and this motivation also explains the choice of studying these specific statistics). 
These results are omitted in the present article. 
On the other hand, we give below some new results: 
namely, we now provide the asymptotic distribution of our statistics 
in the regime with fixed $\theta$, which are Gaussian in three cases, and yield a beta distribution in the fourth case.

In addition, the method used to establish the first statement on each statistic, 
giving the probability of an elementary event, is different from~\cite{AuBoNiPi2016}. 
While we relied in~\cite{AuBoNiPi2016} on a decomposition of the words representing permutations and on some technical lemmas, 
we now derive these results easily from the generative processes of \cref{sec:process_diagrams,sec:process_words} for \recbiased permutations. 

In the sequel, we shall use the properties of the so-called \emph{digamma} function~\cite{OlverLozierBoisvertClark2010}.
The digamma function is defined by $\Psi(x)=\Gamma'(x)/\Gamma(x)$, with $\Gamma$ the classical gamma function generalizing the factorial. It satisfies the recurrence $\Psi(x + 1) = \Psi(x) +\frac{1}{x}$ which leads to the identity
\begin{equation}\label{eq:Psi}
    \sum_{i=0}^{n-1}\frac1{x + i} = \Psi(x + n)-\Psi(x).
\end{equation} 
As $x\to +\infty$, it admits the asymptotic expansion  $\Psi(x) = \log(x) - \frac{1}{2x} -\frac{1}{12x^2}  +  o\left(\frac{1}{x^2}\right)$. 

%%%%%%%%%%%%%%%%%%%%%%%%%%%%%%%%%%%%%%%%
\subsection{Four statistics and their expectations} \label{subsec:stat_comme_avant}

Throughout this section, we denote by $\EE_n[\cdot]$ the expectation of a random variable on \recbiased permutations of size $n$ for the parameter $\theta$.

\subsubsection*{Number of records}

We start  by computing how the value of parameter $\theta$ influences the expected number of records.

\begin{theorem}\label{lem:records}
Among \recbiased permutations of size $n$ for the parameter $\theta$, 
for any $i \in [n]$, 
the probability that there is a record at position $i$ is: 
$\PP_{n}(\text{record at }i)  = \frac{\theta}{\theta  +  i - 1}$.
\end{theorem}

\begin{proof}
There is a record at $i$ if and only if the point inserted in the $i$-th column by the process of \cref{sec:process_diagrams} 
is the highest among those already inserted, and this happens with probability $\frac{\theta}{\theta  +  i - 1}$. 
\end{proof}

\begin{corollary} \label{thm:nb_of_records}
Among \recbiased permutations of size $n$ for the parameter $\theta$, the expected value of the number of records is: 
$\EE_{n}[\rec] = \sum_{i=1}^{n}\frac{\theta}{\theta + i-1} = \theta ( \Psi(\theta + n)-\Psi(\theta))$. 
\end{corollary}

\begin{proof}
The expression of $\EE_{n}[\rec]$ is simply obtained summing over $i$ the result of \cref{lem:records}, and applying \cref{eq:Psi} for the right hand term. 
\end{proof}

\begin{remark} \label{rk:cycles_records}
\cref{thm:nb_of_records} can also be recovered from known results on the Ewens distribution. 
Indeed, as we have seen in \cref{sec:Foata}, the \recbiased distribution on $\sym_n$ is the image, by the fundamental bijection $\bij$ mapping cycles to records, of the Ewens distribution on $\sym_n$. 
Consequently, \cref{thm:nb_of_records} is just a consequence of the well-known expectation 
of the number of cycles under the Ewens distribution (see for instance~\cite[\S 5.2]{Arratia}).  % C'est le paragraphe intitulé ``number of cycles'' 
\end{remark}

\subsubsection*{Number of descents}

Recall that a permutation $\sigma$ of $\sym_{n}$ has a \emph{descent} at position $i\in\{2,\ldots,n\}$ if
$\sigma(i-1)>\sigma(i)$. We denote by $\desc(\sigma)$ the number of descents in $\sigma$. 

\begin{theorem}\label{lem:descent i}
Among \recbiased permutations of size $n$ for the parameter $\theta$, 
for any $i\in\{2,\ldots,n\}$, 
the probability that there is a descent at position $i$ is: 
$\PP_{n}\big(\sigma(i-1)>\sigma(i)\big) =   \frac{(i-1)(2\theta  +  i-2)}{2(\theta + i-1)(\theta + i-2)}$. 
\end{theorem}

\begin{proof}
For a permutation $\sigma$ of size $n$ and any $a \leq n$, we denote by $\norm(\sigma(a))$ the rank of $\sigma(a)$ in $\sigma(1) \sigma(2) \dots \sigma(a)$, 
that is to say $\norm(\sigma(a)) = j$ exactly when $\sigma(a)$ is the $j$-th smallest element in $\{\sigma(1), \sigma(2), \dots, \sigma(a)\}$. 
Splitting according to the possible values of $\norm(\sigma(i-1))$ and $\norm(\sigma(i))$, we have
\[
\PP_{n}\big(\sigma(i-1)>\sigma(i)\big) = \sum_{j=1}^{i-1} \sum_{k=1}^{j} \PP_n\big(\norm(\sigma(i-1))=j \text{ and } \norm(\sigma(i)=k)\big).
\]
From the process of \cref{sec:process_diagrams} for generating permutations as diagrams, 
we know that the rank of $\sigma(i-1)$ is independent from the rank of $\sigma(i)$, 
and we know the probabilities of these ranks taking any given value. 
We therefore obtain
\begin{align*}
\PP_{n}\big(\sigma(i-1)>\sigma(i)\big) & = \sum_{j=1}^{i-1} \sum_{k=1}^{j} \PP_n\big(\norm(\sigma(i-1))=j\big) \cdot \PP_n\big(\norm(\sigma(i)=k)\big) \\
& = \sum_{j=1}^{i-1} \sum_{k=1}^{j} \PP_n\big(\norm(\sigma(i-1))=j\big) \cdot \frac{1}{\theta+i-1} \\
& = \Big(\sum_{j=1}^{i-2} \sum_{k=1}^{j} \frac{1}{\theta+i-2} + \sum_{k=1}^{i-1} \frac{\theta}{\theta+i-2}\Big)\cdot \frac{1}{\theta+i-1} \\
& = \frac{1}{(\theta+i-2) (\theta+i-1)} \big( (\sum_{j=1}^{i-2} j) + \theta (i-1) \big) \\
& = \frac{1}{(\theta+i-2) (\theta+i-1)} \big( \tfrac{(i-1)(i-2)}{2}+ \theta (i-1) \big) \\
& = \frac{(i-1)(2\theta  +  i-2)}{2(\theta + i-1)(\theta + i-2)},
\end{align*}
concluding the proof. 
\end{proof}

\begin{corollary}\label{thm:nb_of_descents}
Among \recbiased permutations of size $n$ for the parameter $\theta$, 
the expected value of the number of descents is: 
$\EE_{n}[\desc]= \frac{n(n-1)}{2(\theta  + n-1)}$.
\end{corollary}

\begin{proof}
The partial fraction decomposition of the formula for $\PP_{n}\big(\sigma(i-1)>\sigma(i)\big)$ in \cref{lem:descent i} gives:
% \[
$2\PP_{n}(\sigma(i-1)>\sigma(i)) = 1 +  \frac{\theta(\theta-1)}{\theta + i-1} - \frac{\theta(\theta-1)}{\theta + i-2}$.
% \]
Hence, we have a telescopic series when summing for $i$ from $2$ to $n$, which yields:
% \[
$2\EE_{n}[\desc] = (n-1)  +  \frac{\theta(\theta-1)}{\theta + n-1} - \frac{\theta(\theta-1)}{\theta}$.
% \]
This gives the announced expression for $\EE_{n}[\desc]$ after some elementary simplifications. 
%The asymptotic behavior (for $n \to \infty$ with $\theta$ fixed) follows immediately from this expression.
\end{proof} 

\begin{remark}\label{rk:descents_exc}
Similarly to \cref{rk:cycles_records}, the number of descents among \recbiased permutations of size $n$ 
is mapped by the fundamental bijection $\bij^{-1}$ to a statistics on Ewens permutations of the same size, 
specifically, to the number of anti-excedances. 
An anti-excedance of $\sigma \in \sym_n$ is an index $i \in [n]$ such that $\sigma(i) <i$. 
The proof that descents of $\sigma$ are equinumerous with anti-excedances of $\bij^{-1}(\sigma)$ is a simple adaptation of the proof of Theorem 1.36 in~\cite{Bona}, p. 110--111 
(and is omitted here).
To our knowledge, this statistics was not studied on Ewens permutations. However, a related statistics was, 
namely, the number of weak excedances -- see~\cite[Thm 1.2 and \S 5]{Valentin}. A weak excedance of $\sigma \in \sym_n$ is an index $i \in [n]$ such that $\sigma(i) \geq i$. 
Clearly, the numbers of anti-excedances and of weak excedances in a permutation of size $n$ sum up to $n$, 
and results on one of these statistics yields results on the other. 
In particular, the expression of $\EE_{n}[\desc]$ in \cref{thm:nb_of_descents} easily follows from the first formula in~\cite[\S 5.4]{Valentin} 
(summing over $i$ and taking the difference with $n$) or from the first equation of~\cite[Thm. 1.2]{Valentin} .   
\end{remark}

\subsubsection*{Number of inversions}

Recall that an \emph{inversion} in a permutation $\sigma\in\sym_{n}$ is a pair $(i,j)\in[n]\times[n]$ such that $i<j$ and $\sigma(i) > \sigma(j)$. 
In the word representation of permutations, this corresponds to a pair of elements in which the largest is to the left of the smallest. 
For any $\sigma\in\sym_{n}$, we denote by $\inv(\sigma)$ the number of inversions of $\sigma$, 
and by $\inv_j(\sigma)$ the number inversions of the form $(i,j)$ in $\sigma$, for any $j\in [n]$. 
More formally, $\inv_{j}(\sigma) =\big| \{i\in[j-1]:\ (i,j)\text{ is an inversion of }\sigma\}\big|$.

\begin{theorem}\label{lem:inversion_at_j}
Among \recbiased permutations of size $n$ for the parameter $\theta$, 
for any $j \in [n]$ and $k\in[0,j-1]$, 
the probability that there are $k$ inversions of the form $(i,j)$ is: \\
$\PP_{n}\big(\inv_j(\sigma) =k\big) =  \frac{1}{\theta+j-1}$ if $k\neq 0$ 
and $\PP_{n}\big(\inv_j(\sigma) =k\big) =\frac{\theta}{\theta+j-1}$ if $k=0$. 
\end{theorem}

\begin{proof}
By definition, the value of $\inv_j(\sigma)$ depends only on how $\sigma(j)$ compares to the $\sigma(i)$ for $i<j$. 
More precisely, $\inv_j(\sigma) = k$ if and only if $\sigma(j)$ is the $(j-k)$-th largest element of $\sigma$ among the first $j$. 
From \cref{sec:process_diagrams}, the probability that $\sigma(j)$ is the largest element of $\sigma$ among the first $j$ is $\frac{\theta}{\theta+j-1}$, 
and this proves the statement of the theorem in the case $k=0$. 
For any value $k\neq 0$ with $k <j$, from \cref{sec:process_diagrams}, the probability that $\sigma(j)$ is the $(j-k)$-th largest element of $\sigma$ among the first $j$ 
is $\frac{1}{\theta+j-1}$, concluding the proof.
\end{proof}

\begin{corollary}\label{thm:nb_of_inversions}
Among \recbiased permutations of size $n$ for the parameter $\theta$, 
the expected value of the number of inversions is: 
$\EE_{n}[\inv] = \frac{n(n + 1-2\theta)}4 + \frac{\theta(\theta-1)}{2}\left(\Psi(\theta + n)-\Psi(\theta)\right)$. \\
\end{corollary}

\begin{proof}
We first compute $\EE_{n}[\inv_{j}]$ for any $j \in [n]$. With \cref{lem:inversion_at_j} and noting that the maximum possible value of $\inv_j$ is $j-1$, we have 
\begin{equation} \label{eq:esp_inv}
\EE_{n}[\inv_{j}] 
= \sum_{k=0}^{j-1}k\cdot\PP_{n}(\inv_{j}=k)
= \sum_{k=1}^{j-1}k\cdot\frac{1}{\theta+j-1} = 
\frac{j(j-1)}{2(\theta  +  j - 1)}.
\end{equation}
Consequently, $\EE_{n}[\inv] = \sum_{j=1}^{n}\EE_{n}[\inv_{j}] = \sum_{j=1}^{n} \frac{j(j-1)}{2(\theta  +  j - 1)}$. 
Using the properties of the digamma function $\Psi$ reviewed at the beginning of \cref{sec:stat}, 
and noting that $\frac{j(j-1)}{\theta+j-1} = j- \theta + \frac{\theta(\theta-1)}{\theta+j-1}$, 
this sum may be expressed as 
$
\EE_{n}[\inv] = \frac{n(n + 1-2\theta)}4 + \frac{\theta(\theta-1)}{2}\left(\Psi(\theta + n)-\Psi(\theta)\right)
$.
\end{proof}		

\begin{remark}
Unlike for the previous statistics, we are not aware of any natural interpretation of the number of inversions of $\sigma$ in $\bij^{-1}(\sigma)$. 
Therefore, there is no known equivalent of the above statements for the Ewens distribution. 
\end{remark}

\subsubsection*{First value}

Naturally, the behavior of the first value $\sigma(1)$ in \recbiased permutations 
differs completely from the three statistics studied above. (If needed, see \cref{eq:rising} for the definition of the rising factorial $\rfact{\theta}{n}$.)

\begin{theorem}\label{cor:first_value=k}
Among \recbiased permutations of size $n$ for the parameter $\theta$, 
for any $k \in [n]$, 
the probability that a permutation starts with $k$ is: 
$\PP_{n}(\sigma(1)=k) = \frac{(n-1)!\,\rfact{\theta}{n-k}\theta}{(n-k)!\rfact{\theta}{n}}$.
\end{theorem}

\begin{proof}
The statement is proved using the generative process of \cref{sec:process_words} (\cref{fig:array-tree}) for permutations viewed as words. 
In this process, for the first element to be equal to $k$, it means that the first $k-1$ insertions (of the values $1$ to $k-1$) correspond to non-records, 
and that the $k$-th insertion (of the value $k$) corresponds to a record. This happens with the following probability: 
\[
\PP_{n}(\sigma(1)=k) =  \prod_{i=1}^{k-1} \frac{n-i}{\theta+n-i} \cdot \frac{\theta}{\theta+n-k} = \frac{(n-1)!\,\rfact{\theta}{n-k}\theta}{(n-k)!\rfact{\theta}{n}}. \qedhere
\]
\end{proof}		

\begin{remark}
The above theorem was also recently proved by B. Corsini in~\cite[Corollary 2.3]{Corsini}. Both our and his proof rely on the generative process of \cref{sec:process_words}.
\end{remark}

\begin{corollary}\label{thm:first_value}
Among \recbiased permutations of size $n$ for the parameter $\theta$, 
the expected value of the first element of a permutation is: 
$\EE_{n}[\sigma(1)] = \frac{\theta + n}{\theta + 1}$.
\end{corollary}

\begin{proof}
From~\cref{cor:first_value=k}, we have 
\[
 \EE_{n}[\sigma(1)] = \sum_{k=1}^{n}k\,\PP_{n}(\sigma(1)=k) = \frac{(n-1)!\theta}{\rfact{\theta}{n}} \sum_{k=1}^{n}k \frac{\rfact{\theta}{n-k}}{(n-k)!} 
 = \frac{(n-1)!\theta}{\rfact{\theta}{n}} \, \sum_{j=0}^{n-1} (n-j) \frac{\rfact{\theta}{j}}{j!}.
\]
By the binomial formula, $[z^{m}](1-z)^{-\alpha}=\frac{\rfact\alpha m}{m!}$, where as usual $[z^m]F(z)$ denotes the coefficient of $z^m$ in the series expansion of $F(z)$. 
Observe also that $(n-j) = [z^{n-j-1}] (1-z)^{-2}$. 
Therefore 
\begin{align*}
\EE_{n}[\sigma(1)] & = \frac{(n-1)!\theta}{\rfact{\theta}{n}} \, \sum_{j=0}^{n-1} [z^{n-j-1}] (1-z)^{-2} \cdot[z^{j}](1-z)^{-\theta} \\ 
& =  \frac{(n-1)!\theta}{\rfact{\theta}{n}} \, [z^{n-1}](1-z)^{-\theta-2} = \frac{(n-1)!\theta}{\rfact{\theta}{n}} \, \frac{\rfact{(\theta+2)} {n-1}}{(n-1)!},
\end{align*}
which yields the announced result after elementary simplifications.
\end{proof}

\begin{remark}
Our proof of \cref{thm:first_value} relies on calculus, but gives a very simple expression for $\EE_{n}[\sigma(1)]$. 
We could therefore hope for a more combinatorial proof of \cref{thm:first_value}, but we were not able to find it, and leave it as an open question.
\end{remark}

\begin{remark}
Under the inverse of the fundamental bijection, the first value of $\sigma$ is mapped to the minimum, over all cycles, of the maximal value in a cycle of $\bij^{-1}(\sigma)$. 
To the best of our knowledge, the behavior of this statistics on Ewens permutations has not been previously studied. 
Therefore, our results above can also be interpreted as new results on Ewens permutations. 
\end{remark}

%%%%%%%%%%%%%%%%%%%%%%%%%%%%%%%%%%%%%%%%
\subsection{Different regimes for different $\theta$}

In this short section, we compute the asymptotic equivalents of the expectations we just obtained for different regimes of $\theta$: constant, sublinear $\theta=n^\varepsilon$ with $\varepsilon\in(0,1)$, linear $\theta=\lambda n$ for $\lambda>0$ and superlinear $\theta=n^\gamma$ with $\gamma>1$. The computations are all straightforward, using the asymptotic behavior of $\Psi$ when needed. The results are given in the table below.

\begin{table}[ht]\label{table:asymptotics}
\begin{center}
\begin{tabular}{l|l|l|l|l|l|c}
& $\theta=1$ & fixed $\theta>0$ & $\theta = n^{\varepsilon}$,  & $\theta =  \lambda n$, & $\theta = n^{\delta}$, &See \\
& {\small (uniform)}& & {\small $0<\varepsilon<1$ }& {\small $\lambda>0$} & {\small $\delta>1$ }&Cor.\\
\hline 
$\EE_{n}[\rec]$ & $ \log n $ & $ \theta\cdot \log n $ & $ (1-\varepsilon)\cdot n^{\varepsilon}\log n$ & $ \lambda \log(1 + 1/\lambda)\cdot n$ & $ n$ & \ref{thm:nb_of_records}\\
$\EE_{n}[\desc]$ & $  n/2 $ & $  n/2 $ & $  n/2$ & $ n/2(\lambda + 1)$ & $ n^{2-\delta}/2$ & \ref{thm:nb_of_descents}\\
$\EE_{n}[\inv]$ & $ n^2/4$ & $ n^2/4$ & $ {n^{2}}/4$ & $ {n^{2}}/4 \cdot f(\lambda)$ & $ {n^{3-\delta}}/{6}$ & \ref{thm:nb_of_inversions}\\
$\EE_{n}[\sigma(1)]$ & $ n/2$ & $  n/(\theta  + 1) $ & $  n^{1-\varepsilon}$ & $ (\lambda + 1)/\lambda$ & $ 1$ & \ref{thm:first_value}\\
\end{tabular}
\end{center}
\caption{Asymptotic behavior in expectation of some permutation statistics on record-biased permutations for the parameter $\theta$.  
We use the shorthand $f(\lambda) = 1-2\lambda + 2\lambda^{2}\log\left(1 + 1/\lambda\right)$.
All the results in this table are asymptotic equivalents.
}
\end{table}

%%%%%%%%%%%%%%%%%%%%%%%%%%%%%%%%%%%%%%%%
\subsection{Limit laws for fixed $\theta$}

In this section, and this section only, we consider the case where $\theta>0$ is fixed. Under this classical assumption, we establish that the number of records, the number of descents and the number of inversions are asymptotically Gaussian, whereas the first value tends to a beta distribution of parameters $(1,\theta)$.
For the first two results, we rely on  known properties of the Ewens distribution.

\begin{figure}[ht]
\centering
\begin{minipage}{.47\textwidth}
  \centering
  \includegraphics[trim={0 0 0 40pt}, clip, width=\linewidth]{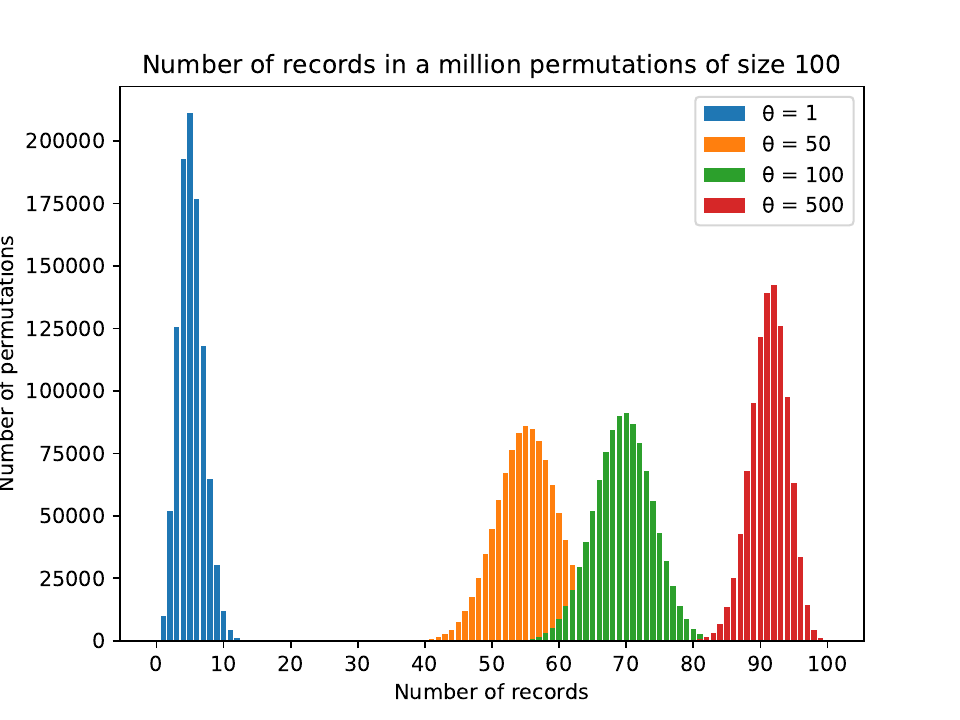}
  \captionof{figure}{Empirical distribution of the number of records in~$10^6$ \recbiased permutations of size~$100$ for several values of~$\theta$.}
  \label{fig:histo_records}
\end{minipage}\hfill
\begin{minipage}{.47\textwidth}
  \centering
  \includegraphics[trim={0 0 0 40pt}, clip, width=\linewidth]{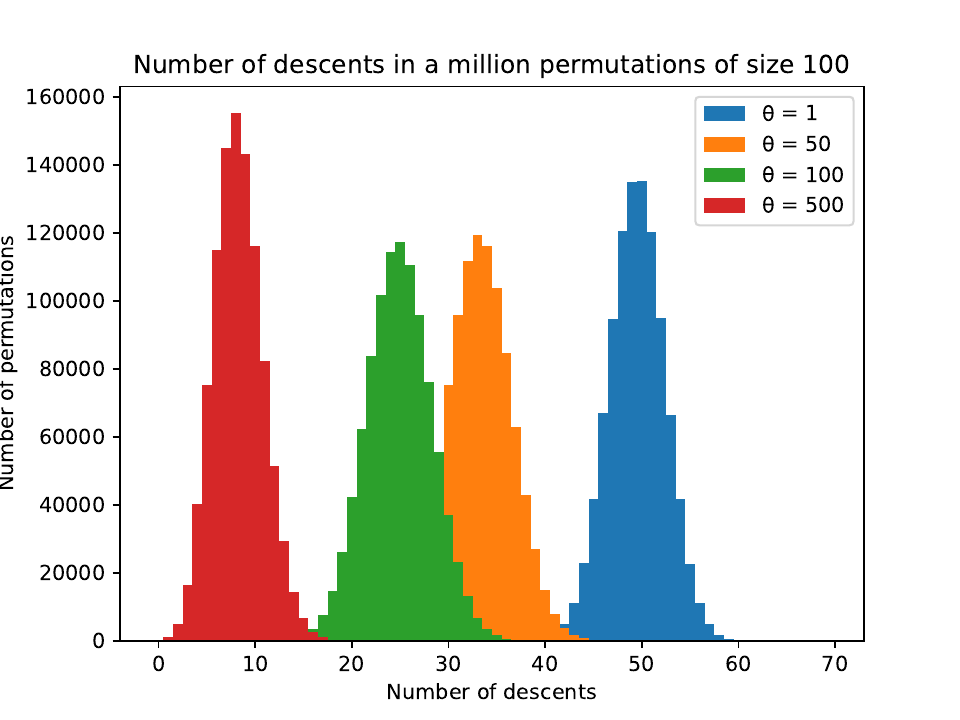}
  \captionof{figure}{Empirical distribution of the number of descents in~$10^6$ \recbiased permutations of size~$100$ for several values of~$\theta$.}
  \label{fig:histo_descents}
\end{minipage}
\end{figure}

\begin{figure}[ht]
\centering
\begin{minipage}{.47\textwidth}
  \centering
  \includegraphics[trim={0 0 0 40pt}, clip, width=\linewidth]{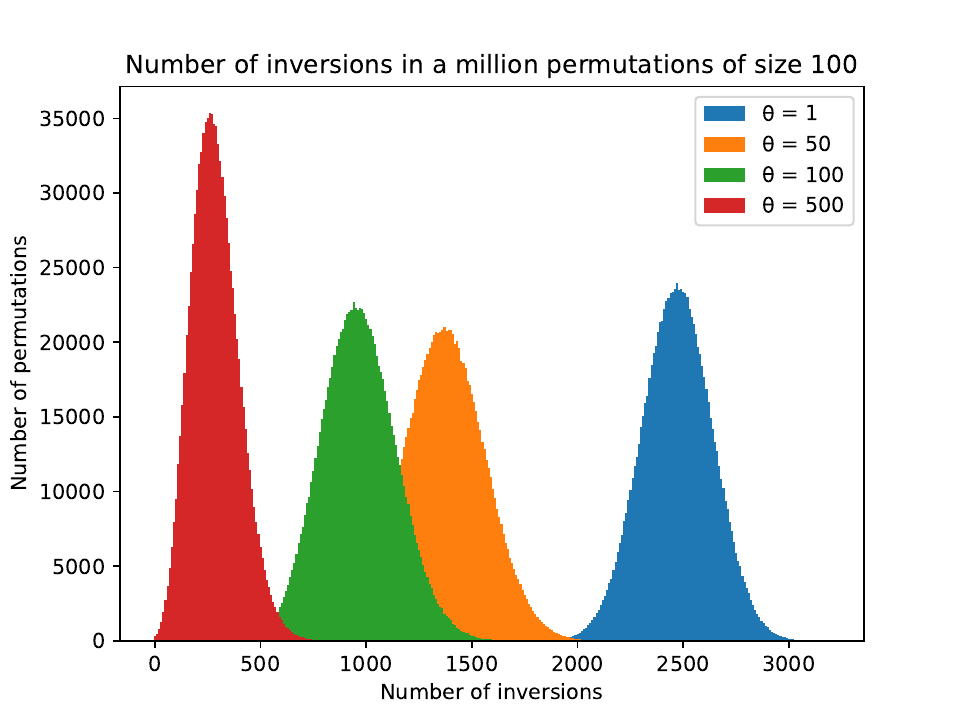}
  \captionof{figure}{Empirical distribution of the number of inversions in~$10^6$ \recbiased permutations of size~$100$ for several values of~$\theta$.}
  \label{fig:histo_inversions}
\end{minipage}\hfill
\begin{minipage}{.47\textwidth}
  \centering
  \includegraphics[trim={0 0 0 40pt}, clip, width=\linewidth]{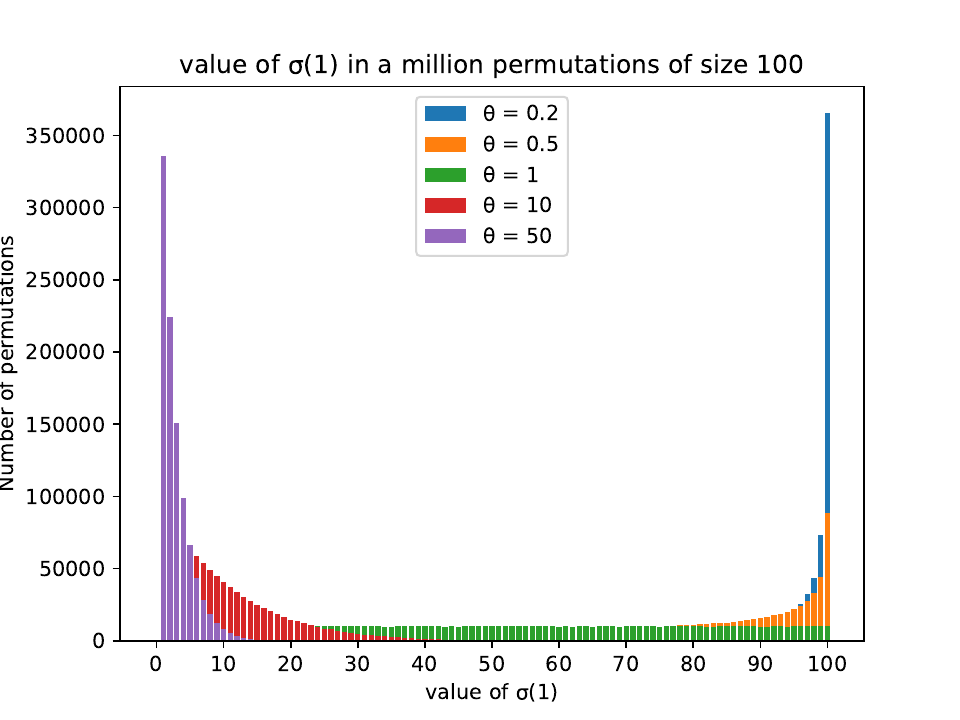}
  \captionof{figure}{Empirical distribution of the value of~$\sigma(1)$ in~$10^6$ \recbiased permutations of size~$100$ for several values of~$\theta$.}
  \label{fig:histo_sigma1}
\end{minipage}
\end{figure}

These convergences appear clearly on the simulations we performed using the random generator described in Section~\ref{sec:random generation}, as can be seen on Fig.~\ref{fig:histo_records} for the number of records, Fig.~\ref{fig:histo_descents} for the number of descents, Fig.~\ref{fig:histo_inversions} for the number of inversions and Fig.~\ref{fig:histo_sigma1} for the first value.

%%%%%%%%%%%%%%%%%%%%%%

\subsubsection*{Number of records}

\begin{theorem}
The number of records in \recbiased permutations is asymptotically normal. More precisely, 
letting $R_n$ be the random variable which denotes the number of records in a \recbiased permutation of size $n$ for any fixed value of the parameter $\theta$, 
we have  
\[
 \frac{R_n - \theta \log(n)}{\sqrt{\theta \log(n)}}  \stackrel{(d)}{\to} \mathcal{N}(0,1) .
\]
\end{theorem}
As usual, the notation $\mathcal{N}(0,1)$ above denotes the normal distribution with mean $0$ and variance~$1$, 
and $ \stackrel{(d)}{\to} $ denotes the convergence in distribution. 

\begin{proof}
As in \cref{rk:cycles_records}, it follows immediately from the distribution of the number of cycles in Ewens permutations, see for instance~\cite[\S 5.2]{Arratia} and more precisely Equation (5.22) p.103. 
\end{proof}

%%%%%%%%%%%%%%%%%%%%%%
\subsubsection*{Number of descents}

\begin{theorem}
The number of descents in \recbiased permutations is asymptotically normal. More precisely, 
letting $D_n$ be the random variable which denotes the number of descents in a \recbiased permutation of size $n$ for any fixed value of the parameter $\theta$, 
we have  
\[
 \frac{D_n - n/2}{\sqrt{n/12}}  \stackrel{(d)}{\to} \mathcal{N}(0,1).
\]
\end{theorem}

\begin{proof}
As explained in \cref{rk:descents_exc}, the distribution of $D_n$ is the same as the distribution of the number of anti-excedances in Ewens permutations, 
and $D_n$ is therefore distributed like $n-W_n$, where $W_n$ denotes the number of weak excedances in a random Ewens permutation of size $n$. 
Since $W_n$ is asymptotically normal (see~\cite[Thm 1.2]{Valentin}), 
it follows that $D_n$ also is. 
To conclude the proof, we just need the asymptotic behaviors of expectation and variance of $D_n$. 
The expectation behaves like $n/2$ as shown in \cref{thm:nb_of_descents}. 
The variance of $D_n$, which is equal to that of $W_n$, can be computed from the results of~\cite{Valentin}. 
More precisely, with the notation of~\cite[Eq (5.3) and immediately after]{Valentin} it behaves asymptotically as $K(1,1) \cdot n$ and elementary computations give $K(1,1) = 1/12$. 
\end{proof}

%%%%%%%%%%%%%%%%%%%%%%
\subsubsection*{Number of inversions}

\begin{theorem}
The number of inversions in \recbiased permutations is asymptotically normal. More precisely, 
letting $I_n$ be the random variable which denotes the number of inversions in a \recbiased permutation of size $n$ for any fixed value of the parameter $\theta$, 
we have  
\[
 \frac{I_n - n^2/4}{\sqrt{n^3/36}}  \stackrel{(d)}{\to} \mathcal{N}(0,1).
\]\end{theorem}

\begin{proof}
Fix a positive real number $\theta$. From \cref{sec:process_diagrams}, we observe that (for $\sigma$ a \recbiased permutation of size $n$ for the parameter $\theta$), 
the random variables $\inv_j(\sigma)$ and $\inv_{j'}(\sigma)$ are independent for $j\neq j'$. 
Therefore $I_n$ is the sum of the independent random variables $\inv_j(\sigma)$ for $j$ from $1$ to $n$, whose distributions are given by \cref{lem:inversion_at_j}.  
This is a classical situation in which asymptotic normality can be proved. 
More precisely, we use here Theorem 27.3 (p. 385) in~\cite{Billingsley} (specifically, we show the condition (27.16) with $\delta =1$). 

To this end, we recall from the proof of \cref{thm:nb_of_inversions} (see \cref{eq:esp_inv}) that $\EE_{n}[\inv_{j}] =\frac{j(j-1)}{2(\theta  +  j - 1)}$, 
which implies $\EE_n[I_n] \sim n^2/4$ as we have seen. 
Similarly, we compute 
\begin{equation*} 
\EE_{n}[\inv_{j}^2] 
= \sum_{k=0}^{j-1}k^2\cdot\PP_{n}(\inv_{j}=k)
= \sum_{k=1}^{j-1}k^2\cdot\frac{1}{\theta+j-1} = 
\frac{j(j-1)(2j-1)}{6(\theta  +  j - 1)}.
\end{equation*}
Denoting $\VV_n(\cdot)$ the variance of a random variable on \recbiased permutations of size $n$, it follows that 
\begin{equation*}
\VV_n(\inv_j) = \frac{j(j-1)(2j-1)}{6(\theta  +  j - 1)} - \frac{j^2(j-1)^2}{4(\theta  +  j - 1)^2} = \frac{j(j-1)(j^2 +(4\theta-3)j +2-2\theta)}{12(\theta+j-1)^2}, 
\end{equation*}
so that $\VV_n(I_n) \sim n^3/36$ as $n \to \infty$.

Using Theorem 27.3 (p. 385) in~\cite{Billingsley} with the condition (27.16) for $\delta =1$, 
to conclude the proof, it is enough to show that 
\[
 \sum_{j=1}^n \frac{1}{\sqrt{\sum_{j=1}^n \VV_n(\inv_j)}^3} \EE_n[\inv_j^3] \to 0 \text{ as } n \to \infty.
\]
With computations similar to those for $\EE_n[\inv_j]$ and $\EE_n[\inv_j^2]$, we see that $\EE_n[\inv_j^3]$ is of order $j^3$, so that $\sum_{j=1}^n \EE_n[\inv_j^3]$ is of order $n^4$. 
On the other hand, $\sqrt{\sum_{j=1}^n \VV_n(\inv_j)}^3 = \VV_n(I_n)^{3/2}$ is of order $n^{9/2}$, 
and the proof is complete. 
\end{proof}

%%%%%%%%%%%%%%%%%%%%%%
\subsubsection*{First value}

\begin{theorem}
Fix $\theta$ a positive real number. For $\sigma$ a \recbiased permutation of size $n$ for the parameter $\theta$, 
the random variable $\sigma(1)/n$ is asymptotically distributed following a beta distribution of parameters $(1,\theta)$.
\end{theorem}

\begin{proof}
For any positive integer $r$, the $r$-th moment of a beta distribution of parameters $(1,\theta)$ is $\frac{r!}{(\theta+1)^{(r)}}$ (see \emph{e.g.} \cite[Ch. 25, p.217]{beta}). 
In particular, the associated moment generating function has positive radius of convergence, so that the beta distribution is determined by its moments (see \emph{e.g.} \cite[Thm 30.1]{Billingsley}).
% see also https://en.wikipedia.org/wiki/Beta_distribution#Higher_moments
Therefore, it is enough to compute the limits of the moments of $\sigma(1)/n$, for~$\sigma$ a \recbiased permutation of size $n$ for the parameter $\theta$, 
and to observe that these limits are indeed $\frac{r!}{(\theta+1)^{(r)}}$.

To establish this, we use two properties of the Eulerian polynomials (which can be found \emph{e.g.} in the recent book~\cite[\S 1.4 and 1.5]{Petersen}). 
For any $r>0$, we denote by $A_r(z) = a_{r,0} + a_{r,1} z + \dots + a_{r,r-1} z^{r-1}$ the Eulerian polynomial of degree $r-1$. 
First, $a_{r,i}$ is the number of permutations of size $r$ with $i$ descents, so that  we have $a_{r,0} + a_{r,1} + \dots + a_{r,r-1} = r!$. 
Second, the following identity (known as the Carlitz identity) holds, for any $r >0$: $\sum_n n^r z^n = \frac{z A_r(z)}{(1-z)^{r+1}}$. 

We can now proceed to the computation of the $r$-th moment of $\sigma(1)$, for $r>0$. (From the first to the second line, we use again the binomial formula as in the proof of \cref{thm:first_value}.)
\begin{align*}
 \EE_n[\sigma(1)^r] &= \sum_{k=1}^n k^r \PP_{n}(\sigma(1)=k) = \frac{(n-1)!\theta}{\rfact{\theta}{n}} \sum_{k=1}^n k^r \frac{\rfact{\theta}{n-k}}{(n-k)!} 
  = \frac{(n-1)!\theta}{\rfact{\theta}{n}} \sum_{j=0}^{n-1} (n-j)^r \frac{\rfact{\theta}{j}}{j!}  \\
 & = \frac{(n-1)!\theta}{\rfact{\theta}{n}} \sum_{j=0}^{n-1} [z^{n-j}] \frac{z A_r(z)}{(1-z)^{r+1}}\cdot [z^j] \frac{1}{(1-z)^\theta} \\
 & = \frac{(n-1)!\theta}{\rfact{\theta}{n}} \sum_{j=0}^{n-1} [z^{n-j-1}] \frac{A_r(z)}{(1-z)^{r+1}}\cdot [z^j] \frac{1}{(1-z)^\theta} \\
 & = \frac{(n-1)!\theta}{\rfact{\theta}{n}} [z^{n-1}] \frac{A_r(z)}{(1-z)^{r+1+\theta}} = \frac{(n-1)!\theta}{\rfact{\theta}{n}} [z^{n-1}] \frac{a_{r,0} + a_{r,1} z+ \dots + a_{r,r-1} z^{r-1}}{(1-z)^{r+1+\theta}} \\
 & = \frac{(n-1)!\theta}{\rfact{\theta}{n}} \sum_{j=1}^{r} a_{r,j-1} [z^{n-j}] \frac{1}{(1-z)^{r+1+\theta}}
\end{align*}

Next, using the same trick once more, we compute (for $n$ large enough, namely $n-1 \geq r$) 
\begin{align*}
 \frac{(n-1)!\theta}{\rfact{\theta}{n}}  [z^{n-j}] \frac{1}{(1-z)^{r+1+\theta}} 
 & = \frac{(n-1)!\theta}{\rfact{\theta}{n}}  \frac{(\theta +r +1)^{(n-j)}}{(n-j)!} \\ 
 & = \frac{(n-1) \dots (n-j+1) \, \theta \, (\theta +r +1) \dots (\theta +r+n-j)}{\theta (\theta +1) \dots (\theta +n -1)} \\
 & = \frac{(n-1) \dots (n-j+1) \, (\theta +n) \dots (\theta +r+n-j)}{(\theta +1) \dots (\theta + r)} \\
 & = \frac{n^{j-1} (1 + o(1)) \, n^{r-j+1} (1+o(1))}{(\theta +1)^{(r)}} = \frac{n^r (1+o(1))}{(\theta +1)^{(r)}}.
\end{align*}
Finally, 
\[ 
 \EE_n[\sigma(1)^r] = \sum_{j=1}^{r} a_{r,j-1} \frac{n^r (1+o(1))}{(\theta +1)^{(r)}} \text{ so that } \EE_n[\sigma(1)^r] \sim_{n \to \infty} \frac{r! n^r}{(\theta+1)^{(r)}}. 
\] 
Therefore, dividing by $n^r$, the $r$-th moment of $\sigma(1)/n$ coincides with that of a beta distribution of parameters $(1,\theta)$, concluding the proof. 
\end{proof}

%%%%%%%%%%%%%%%%%%%%%%%%%%%%%%%%%%%%%%%%%%%%%%%%%%%%%%%%%%%%%%%%%%%%%%%%%%%%%%%%
\section{Permuton limit}\label{sec:permuton}

\begin{figure}[ht]
\includegraphics[scale=.18]{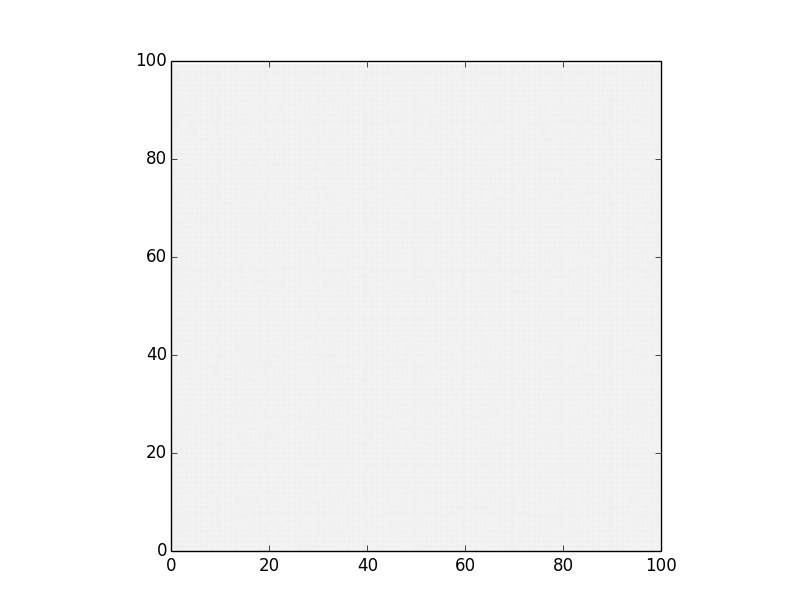}\hspace{-1mm}
\includegraphics[scale=.18]{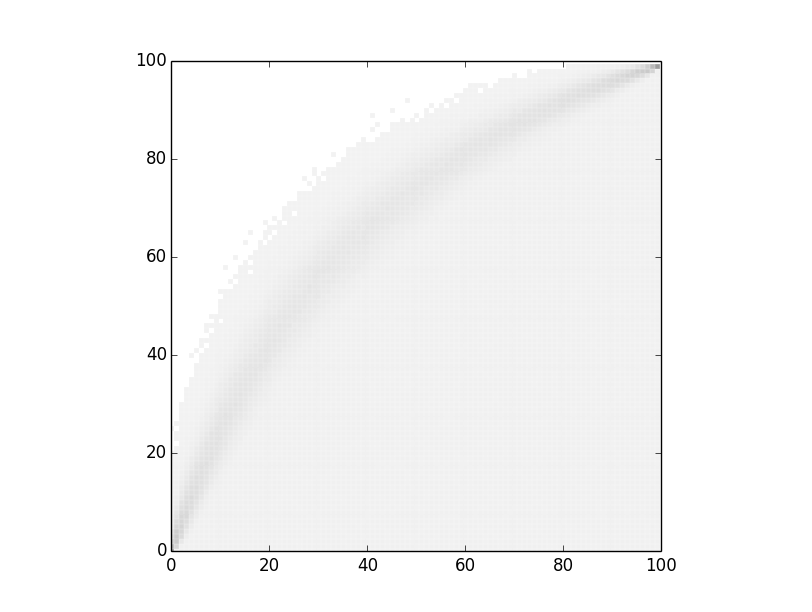}\hspace{-1mm}
\includegraphics[scale=.18]{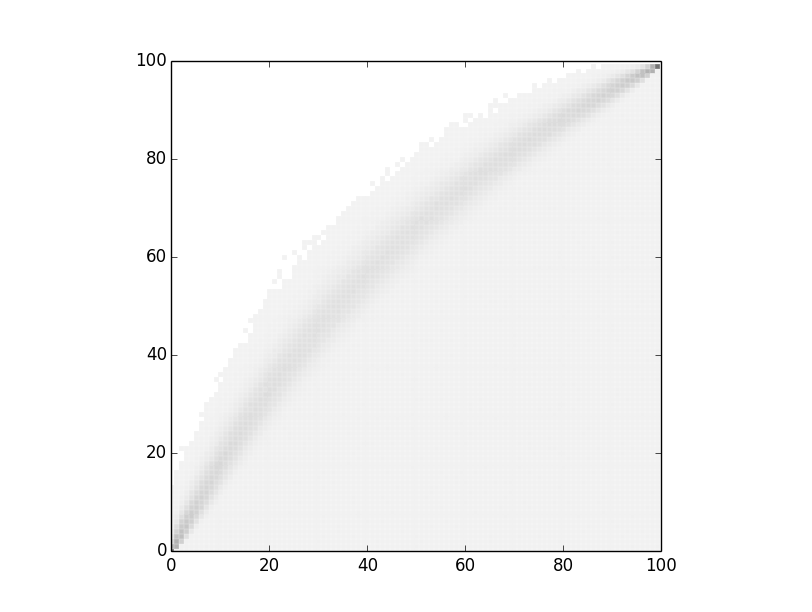}\hspace{-1mm}
\includegraphics[scale=.18]{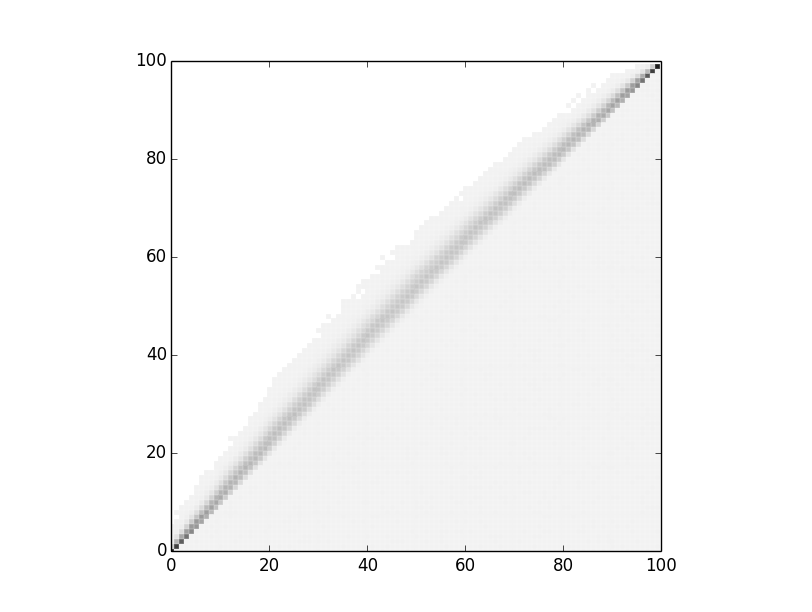}
\caption{Random permutations under the \recbiased distribution on $\sym_{100}$ with, from left to right, $\theta=1$ (corresponding to the uniform distribution)$,50,100,$ and $500$. 
For each diagram,  the darkness of a point $(i,j)$ is proportional to the number of generated permutations $\sigma$ such that $\sigma(i)=j$, for a sampling of 10000 random permutations of size $100$.\label{fig:diagrams}}
\end{figure}

Using our efficient random samplers, we have been able to visualize diagrams of random \recbiased permutations. This is how \cref{fig:diagrams} has been obtained. Explaining the limit shape that seems to emerge from this picture amounts to determining the permuton limit of random \recbiased permutations. This is the purpose of this section.

\subsection{Permutons: definition and criterion of convergence}

The notion of \emph{permutons} was defined independently in~\cite{PresuttiStromquist} and~\cite{Permutons}, 
although the terminology was introduced later in~\cite{FinitelyForcible}. 
It provides a convenient framework to describe the limits of permutations, once viewed as their rescaled diagrams -- see \emph{e.g.} \cite{Winkler,PHC,PHC2}. 

Permutons, denoted $\mu$ here, are probability measures on the unit square $[0,1]^2$ such that 
\begin{equation}\label{eq:marginals}
\text{for all } a \leq b \in [0,1], \quad \mu([a,b]\times[0,1]) = b-a \text{ and } \mu([0,1]\times[a,b]) = b-a. 
\end{equation}
Their projections on the horizontal and vertical axes are the uniform measure on $[0,1]$.
We sometimes say that permutons have \emph{uniform marginals} because of the property of \cref{eq:marginals}.

\begin{figure}[ht]
\begin{center}
 \begin{tikzpicture}[scale=0.3]
\draw (0,0) grid (7,7);
\draw (0.5,1.5) [fill] circle (.2);
\draw (1.5,3.5) [fill] circle (.2);
\draw (2.5,6.5) [fill] circle (.2);
\draw (3.5,0.5) [fill] circle (.2);
\draw (4.5,4.5) [fill] circle (.2);
\draw (5.5,5.5) [fill] circle (.2);
\draw (6.5,2.5) [fill] circle (.2);
\end{tikzpicture}
\qquad  \begin{tikzpicture}[scale=0.3]
\fill[color=lightgray] (0,1) rectangle (1,2);
\fill[color=lightgray]  (1,3) rectangle (2,4);
\fill[color=lightgray]  (2,6) rectangle (3,7);
\fill[color=lightgray]  (3,0) rectangle (4,1);
\fill[color=lightgray]  (4,4) rectangle (5,5);
\fill[color=lightgray]  (5,5) rectangle (6,6);
\fill[color=lightgray]  (6,2) rectangle (7,3);
\draw (0,0) grid (7,7);
\end{tikzpicture}
\end{center}
 \caption{The diagram and the permuton of $\sigma = 2 4 7 1 5 6 3$. \label{fig:permuton_def}}
\end{figure}
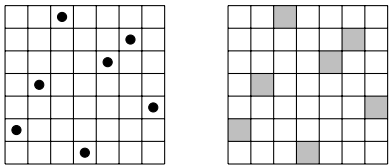

With every permutation $\sigma$, denoting $n$ its size, we can associate a permuton as follows. 
Recall that the \emph{diagram} of $\sigma$ consists of the set of points at coordinates $(i,\sigma(i))$ 
in an $n\times n$ grid. 
Rescaling the diagram to the unit square, and replacing points by small squares carrying total mass $1$, uniformly distributed, 
we obtain the permuton associated with $\sigma$, which we denote $\mu_\sigma$. See \cref{fig:permuton_def}.
More formally, $\mu_\sigma$ is obtained 
dividing the unit square in a regular $n\times n$ grid, 
and distributing uniformly the total mass $1$ on the small squares $[(i-1)/n,i/n]\times[(\sigma(i)-1)/n,\sigma(i)/n]$. 
It is immediate to check that this construction ensures the uniform marginals property, 
so that~$\mu_\sigma$ is indeed a permuton, for every permutation $\sigma$. 

\medskip

Since permutons are measures, the weak convergence topology allows to define the notion of convergence of permutons. 
The distance between permutons is defined by 
\begin{equation}
  d_{all}(\mu , \nu) = \sup (|\mu(R)-\nu(R)|), 
  \label{eq:def_permuton_dist}
\end{equation}
where the supremum is taken over \emph{all} rectangles $R = [a,b]\times[c, d]$ included in the unit square $[0,1]^2$. 
Subsequently, a sequence $(\mu_n)_n$ of permutons converges to a permuton $\mu$ (which we write $\mu_n \to \mu$) when $d_{all}(\mu_n,\mu) \to 0$. 
And for $(\mu_n)_n$ a sequence of \emph{random} permutons, $\mu_n$ converges in probability to $\mu$, which we write $\mu_n \toinp \mu$, when for every $\varepsilon >0$ it holds that $\PP(d_{all}(\mu_n,\mu)>\varepsilon) \to 0$ as $n \to \infty$. 

\medskip

We can define a variant $d(\mu, \nu)$ of the distance between permutons, similarly to \cref{eq:def_permuton_dist}, 
but where we take the supremum only over rectangles $R$ whose bottom-left corner is in $(0,0)$. 
Although the distances differ, it is easy to see that they differ only by a factor of at most $4$ (this is also observed in \cite[eq.(34)]{Permutons}), so that $\mu_n \to \mu$ is equivalently defined by $d(\mu_n,\mu) \to 0$ and $\mu_n \toinp \mu$ is equivalent to having, for every $\varepsilon >0$, $\PP(d(\mu_n,\mu)>\varepsilon) \to 0$ as $n \to \infty$. 

We will be particularly interested in the case where the converging sequence $\mu_n$ is random and obtained from a sequence of random permutations of increasing sizes, 
that is to say, $\mu_n$ is the random permuton associated with a random permutation $\sigma_n$ of size $n$, for any $n$. 
In this case, as the following lemma demonstrates, we can restrict even further the family of rectangles on which we take the supremum. 
Namely, we define
\[
   d_{\cdot /n}(\mu , \nu) = \sup (|\mu(R)-\nu(R)|), 
\]
where $R$ ranges over the family $\mathcal{R}_n$ of rectangles of the form $[0,i/n]\times[0,j/n]$, for $0\leq i,j \leq n$.

\begin{lemma}
For any $n$, let $\sigma_n$ be a random permutation of size $n$, 
and let $\mu_n$ be the associated random permuton. 
The following are equivalent: 
\begin{itemize}
 \item $\mu_n \toinp \mu$; 
 \item for every $\varepsilon \in (0,1/2)$, we have $\PP(d_{\cdot /n} (\mu_n, \mu)>\varepsilon) \to 0$ as $n \to \infty$.
\end{itemize}
In addition, if $(\sigma_n)_n$ is a sequence of (non-random) permutations with $|\sigma_n|=n$ and  $\mu_n$ is the (non-random) associated permuton, we have equivalence between $\mu_n \to \mu$ and $d_{\cdot /n} (\mu_n,\mu) \to 0$ as $n \to \infty$. 
\label{lem:conv_permuton}
\end{lemma}

\begin{proof}
For now, fix some positive integer $n$. Let $R$ be any rectangle within $[0,1]^2$ whose bottom-left corner is in $(0,0)$, and let $\nu$ be any permuton. 
We define $R^\square$ (resp. $R^\blacksquare$) as the largest rectangle contained in $R$ (resp. smallest rectangle which contains $R$) whose sides are supported by the lines of the regular $n\times n$ grid inside the unit square. In other words, we set 
\[
 R^\square=[0,a^\square]\times[0,b^\square] \text{ and } 
 R^\blacksquare = [0,a^\square+\tfrac1n]\times[0,b^\square+\tfrac1n], \text{ for } a^\square = \tfrac1n\lfloor an\rfloor\text{ and }b^\square = \tfrac1n\lfloor bn\rfloor  
\]
(unless $a^\square =a$ (resp. $b^\square=b$), in which case we replace $a^\square+\tfrac1n$ by 
$a^\square$ (resp. $b^\square+\tfrac1n$ by $b^\square$) in the definition of  
$R^\blacksquare$). 

Since $R^\square\subseteq R \subseteq R^\blacksquare$, we have $\nu(R^\square)\leq \nu(R) \leq \nu(R^\blacksquare)$. In addition, 
\[
R^\blacksquare \subseteq R^\square\cup [a^\square,a^\square+\tfrac1n]\times[0,1] \cup
[0,1]\times [b^\square,b^\square+\tfrac1n],
\]
so that $\nu(R^\blacksquare) \leq \nu(R^\square) + \frac2n$. 
In particular, it holds that $|\nu(R) - \nu(R^\square)| \leq \tfrac2n$. 
The second statement of \cref{lem:conv_permuton} then follows readily, since (with notation from this statement) 
\[
 d(\mu_n,\mu) = \sup_{R = [0,a]\times[0,b]}(|\mu_n(R)-\mu(R)|) \leq 
\sup_{R = [0,a]\times[0,b]}(|\mu_n(R^\square)-\mu(R^\square)|+\tfrac4n) = d_{\cdot /n}(\mu_n , \mu) +\tfrac4n.
 \]
 
We now turn to the proof of the first statement (using notation from there), and more precisely to proving that the second item implies the first (the converse being obvious). 

We first introduce notation. Recall that $\mathcal{R}_n$ is the set of rectangles of the form $[0,\tfrac{i}n]\times[0,\tfrac{j}n]$ for integers $i,j\leq n$. 
For $\alpha >0$, let $\MS_n(\alpha)$ and $\MG_n(\alpha)$ be the following properties (recall that~$\mu_n$ in the random permuton associated with the random permutation $\sigma_n$ of size $n$): 
 \begin{align*}
 \MS_n(\alpha) & = \exists R^\square\in\mathcal{R}_n,\ \mu_n(R^\square) < \mu(R^\square) -\alpha \\
 \MG_n(\alpha) & = \exists R^\square\in\mathcal{R}_n,\ \mu_n(R^\square) > \mu(R^\square) + \alpha. 
 \end{align*}
Under the assumption of the second item, for any fixed $\alpha>0$, the properties $\MS_n(\alpha)$ and $\MG_n(\alpha)$ both have a probability tending to $0$ as $n$ tends to $\infty$. Also note that, for any $\beta \leq \alpha$, $\MS_n(\alpha)$ implies $\MS_n(\beta)$, hence $\PP(\MS_n(\alpha)) \leq \PP(\MS_n(\beta))$, and similarly with $\MG_n$. 

Fix some $\varepsilon \in (0,1/2)$. 
We first observe (using the union bound) that 
\begin{multline*}
 \PP\left(\sup_{R = [0,a]\times[0,b]}\left| \mu_n(R)-\mu(R)\right|>\varepsilon\right) =  \PP\left(\exists R = [0,a]\times[0,b],\ |\mu_n(R) - \mu(R)| > \varepsilon \right) \\
\leq \PP\left(\exists R = [0,a]\times[0,b],\ \mu_n(R) - \mu(R) > \varepsilon \right) + \PP\left(\exists R = [0,a]\times[0,b],\ \mu(R) - \mu_n(R) > \varepsilon \right).
\end{multline*}
Next, from the inequalities established above, assuming that 
$\mu_n(R) > \mu(R) + \varepsilon$, we have 
\begin{equation*}
\mu_n(R^\square) +\tfrac2n \geq \mu_n(R) > \mu(R) + \varepsilon \geq \mu(R^\square) + \varepsilon, \text{ hence } \mu_n(R^\square) > \mu(R^\square) + \varepsilon -\tfrac2n.
\end{equation*}
Similarly, if $\mu(R) - \mu_n(R) > \varepsilon$, then $\mu_n(R^\square) < \mu(R^\square) - \varepsilon +\tfrac2n$. As a consequence, 
\[
 \PP\left(\sup_{R = [0,a]\times[0,b]}\left| \mu_n(R)-\mu(R)\right|>\varepsilon\right) \leq \PP\left(\MG_n(\varepsilon -\tfrac2n)\right) + \PP\left(\MS_n(\varepsilon -\tfrac2n)\right).
\]
Therefore, for $n$ large enough (\emph{i.e.}, such that $\tfrac{\varepsilon}2 < \varepsilon -\tfrac2n$), 
\[
 \PP\left(\sup_{R = [0,a]\times[0,b]}\left| \mu_n(R)-\mu(R)\right|>\varepsilon\right) \leq \PP\left(\MG_n(\tfrac{\varepsilon}2)\right) + \PP\left(\MS_n(\tfrac{\varepsilon}2)\right),
\]
ensuring that 
$\PP\left(\sup\limits_{R = [0,a]\times[0,b]}\left| \mu_n(R)-\mu(R)\right|>\varepsilon\right)$ 
tends to $0$ as $n$ tends to $\infty$.
\end{proof}

Finally, we point out an important result about permuton convergence: 
it is equivalent to the convergence of all \emph{pattern densities} -- see \cite{PresuttiStromquist,Permutons} in the deterministic case, or \cite[\S 2]{PHC} for convergence of sequences of \emph{random} permutons. 
In the present work, we do not use this characterization, and refer the interested reader to the above papers for details. 

\subsection{The permuton limit of \recbiased permutations for $\theta = \lambda n$} \label{sec:def_mu}
In this section, we define a permuton $\mu$, which we later prove  (see \cref{thm:permutonLimit} and \cref{ssec:proof_permuton_cv}) to be the limit of the sequence of random permutations $(\sigma_n)_n$ defined as follows: 
for each $n$, $\sigma_n$ is a random record-biased permutation of size $n$, for the parameter $\theta = \lambda n$. 
(The reader may notice that this is a different regime than the main one studied in \cref{sec:stat}, where $\theta$ was taken as a fixed positive real number.)

Let $\lambda >0$ be a parameter. 
Let $f_\lambda: [0,1] \to [0,1]$ be the function defined by $f_\lambda(x) = \frac{x(\lambda+1)}{\lambda +x}$. It is obvious but useful to observe that $f_\lambda$ is continuous and increasing from $f_\lambda(0)=0$ to $f_\lambda(1)=1$. It is therefore invertible, and its inverse is given by $f_\lambda^{-1}(x) =\frac{\lambda x}{1+\lambda-x}$. Similarly, an obvious but useful fact is that $f_\lambda(x) \geq x$ for all $x \in [0,1]$. 

We define $Leb_c$ as the measure on the unit square whose support is the curve of equation $y = f_\lambda(x)$, 
with the mass $1$ distributed uniformly along the $x$-axis: in other words, for all $a$, $Leb_c(\{(x,f_\lambda(x)): x \in [0,a]\})$ has measure $a$. 
Then, we define $\mu_c$ as the measure with density $\frac{\lambda}{\lambda+x}$ with respect to $Leb_c$. 

We also define $\mu_u$ as the measure on the unit square whose support is the area below the curve of equation $y = f_\lambda(x)$, 
and with mass $\frac{1}{\lambda+1} \int_0^1 f_\lambda(x)dx$ uniformly distributed on this support. 

\begin{figure}[t]
\centering   
\begin{tikzpicture}[scale=3.2]
\draw (0,0) rectangle (1,1);
\fill [gray!30, domain=0:1, variable=\x](0, 0) -- plot ({\x}, {\x*1.05/(\x+0.05)}) -- (1, 1) -- (1, 0);
\draw [gray, domain=0:1, variable=\x] plot ({\x}, {\x*1.05/(\x+0.05)});
\end{tikzpicture}
\qquad \qquad    
\begin{tikzpicture}[scale=3.2]
\draw (0,0) rectangle (1,1);
\fill [gray!30, domain=0:1, variable=\x](0, 0) -- plot ({\x}, {\x*1.2/(\x+0.2)}) -- (1, 1) -- (1, 0);
\draw [gray, domain=0:1, variable=\x] plot ({\x}, {\x*1.2/(\x+0.2)});
\end{tikzpicture}
\qquad \qquad 
\begin{tikzpicture}[scale=3.2]
\draw (0,0) rectangle (1,1);
\fill [gray!30, domain=0:1, variable=\x](0, 0) -- plot ({\x}, {\x*2/(\x+1)}) -- (1, 1) -- (1, 0);
\draw [gray, domain=0:1, variable=\x] plot ({\x}, {\x*2/(\x+1)});
\end{tikzpicture}
\caption{From left to right: graphical representation of~$\mu$ for~$\lambda=0.05$, $0.2$ and~$1$.}\label{fig:3_mu}
\end{figure}
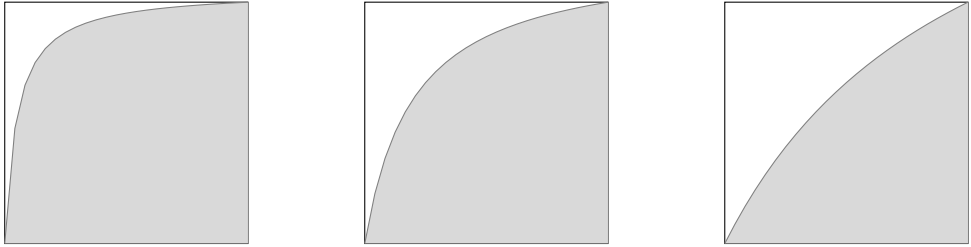

Finally, we let $\mu = \mu_c + \mu_u$ (examples are displayed in \cref{fig:3_mu}).

\begin{remark}
 We mention that the article~\cite{runsort} defines another permuton (called the \emph{runsort permuton}) 
 as the sum of a singular part supported by a curve and an absolutely continuous part supported by the part of the unit square above this curve. 
 Our permuton $\mu$ (up to symmetry) shares this structure with the runsort permuton. 
\end{remark} 

\begin{lemma}
The measure $\mu$ on the unit square defined above is a permuton. 
\end{lemma}

\begin{proof}
We claim that, for each $0\leq a\leq b \leq 1$, $\mu([a,b]\times[0,1]) =b-a$. 
This ensures (for $a=0,b=1$) that $\mu$ is a probability measure on the unit square; 
and it also ensures one part of the uniform marginals property. This claim is easily proved; indeed, we have 
\begin{align*}
 \mu([a,b]\times[0,1]) &= \int_a^b \frac{\lambda}{\lambda+x} dx + \frac{1}{\lambda+1} \int_a^b f_\lambda(x)dx
 = \int_a^b \left(\frac{\lambda}{\lambda+x} +\frac{f_\lambda(x)}{\lambda+1}\right) dx\\ 
 & = \int_a^b \left(\frac{\lambda}{\lambda+x} +\frac{x}{\lambda+x}\right) dx \text{ \quad substituting $f_\lambda(x)$ by its definition} \\
 & = b-a.
\end{align*}

It remains to show that for each $0\leq a\leq b \leq 1$, $\mu([0,1]\times[a,b]) =b-a$. 
For this, we prove below that $\mu([0,1]\times[a,b]) = \mu([1-b,1-a]\times[0,1])$. 
From there, we can use the above computation to derive $\mu([1-b,1-a]\times[0,1])= (1-a) - (1-b) = b-a$. 

To prove $\mu([0,1]\times[a,b]) = \mu([1-b,1-a]\times[0,1])$, 
we shall prove that this identity holds for $\mu_u$ and for $\mu_c$, 
and then deduce it for $\mu = \mu_u + \mu_c$. 

We first observe that the curve of equation $y=f_\lambda(x)$ is symmetric along the line of equation $y=1-x$ (the NorthWest-SouthEast diagonal of the unit square): 
this follows from the equality $f_\lambda(1-f_\lambda(x)) = 1-x$, which is just proved by a simple computation. 
And this is enough to ensure that $\mu_u([0,1]\times[a,b]) = \mu_u([1-b,1-a]\times[0,1])$. 

Now, recalling that $\mu_c$ is the measure with density $\frac{\lambda}{\lambda+x}$ with respect to $Leb_c$, we have 
\begin{align*}
 \mu_c([1-b,1-a]\times[0,1]) & = \int_{1-b}^{1-a} \frac{\lambda}{\lambda+x} dx = 
 \left[ \lambda \log(\lambda+x)\right]_{1-b}^{1-a} = \lambda \log\left(\frac{\lambda+1-a}{\lambda+1-b}\right), \text{ and}\\
 \mu_c([0,1]\times[a,b]) & = \int_{f_\lambda^{-1}(a)}^{f_\lambda^{-1}(b)} \frac{\lambda}{\lambda+x} dx = 
 \left[ \lambda \log(\lambda+x)\right]_{f_\lambda^{-1}(a)}^{f_\lambda^{-1}(b)} = \lambda \log\left(\frac{\lambda+1-a}{\lambda+1-b}\right), 
\end{align*}
where we used that the inverse of $f_\lambda$ is given by $f_\lambda^{-1}(x) =\frac{\lambda x}{1+\lambda-x}$.
This ensures, as wanted, that $\mu_c([0,1]\times[a,b]) = \mu_c([1-b,1-a]\times[0,1])$ and therefore concludes the proof. 
\end{proof}

\begin{remark}
We point out that the permuton $\mu$ (or rather its symmetry w.r.t. a horizontal axis, which we denote $\mu^R$ below) appeared in a similar context in the PhD thesis of J. Borga~\cite[Ch.~5]{Jacopo}, however in disguise. 
\end{remark}

\medskip

We first explain how to regonize $\mu^R$ in~\cite[\S5.2.3]{Jacopo}. By symmetry, from our definition of $\mu$, the permuton $\mu^R$ consists of two parts: 
\begin{itemize}
 \item a first part supported by the curve of equation $y = 1-f_\lambda(x) = 1- \frac{x(\lambda+1)}{\lambda +x}$, of total mass 
 $\int_a^b \frac{\lambda}{\lambda+x} dx = \lambda \log(1+1/\lambda)$; 
 \item and a second part, above this curve, where the mass 
 $ \frac{1}{\lambda+1} \int_0^1 f_\lambda(x)dx = 1 - \lambda \log(1+1/\lambda)$ is uniformly distributed. 
\end{itemize}
On the other hand, the permuton $\nu$ defined in~\cite[\S5.2.3]{Jacopo} (where it is denoted $\nu^\alpha_*$) is described as follows. 
Given a positive $\alpha$, let $\gamma$ be the solution of $\tfrac{1}{\gamma-1}\log(\gamma) = \tfrac{1}{\alpha+1}$. Then $\nu$ consists of two parts: \begin{itemize}
 \item a first part supported by the curve parametrized by $(\tfrac{\gamma^t-1}{\gamma-1},\tfrac{\gamma^{1-t}-1}{\gamma-1})$ for $t \in [0,1]$,
 of total mass $\tfrac{1}{\alpha+1}$; 
 \item and a second part, above this curve, where the mass 
 $\tfrac{\alpha}{\alpha+1}$ is uniformly distributed. 
\end{itemize}
It is straightforward to check that $\mu^R = \nu$, if we take $\alpha = \tfrac{1}{\lambda \log(1+1/\lambda)} -1$, or equivalently $\lambda =\tfrac{1}{\gamma-1}$. 

\medskip

We now describe the context in which $\nu$ appeared in~\cite[Ch.~5]{Jacopo}. 
Specifically, Conjecture 5.2.6 in~\cite{Jacopo} states that $\nu = \mu^R$ is the permuton limit of uniform permutations of size $m=n+k$ with exactly $n$ left-to-right minima and $k$ elements which are not left-to-right minima, with $n$ and $k$ both tending to $\infty$ with $k/n \to \alpha$ for some positive $\alpha$. This is equivalent to stating that $\mu$ is the permuton limit of uniform permutations of size $m=n+k$ with exactly $n$ records and $k$ non-records, in the same regime for $n$ and $k$. 

While this is different from our model of \recbiased permutations, we note that the number of records in the model of~\cite{Jacopo} is in the limit $\frac{1}{\alpha+1}$ times the total size of the permutation. 
In our model, we have seen (see \cref{table:asymptotics}) that the number of records is in the limit $\lambda \log(1+1/\lambda)$ times the total size. 
We observe that, with the relation $\alpha = \tfrac{1}{\lambda \log(1+1/\lambda)} -1$ given earlier, $\frac{1}{\alpha+1} = \lambda \log(1+1/\lambda)$. 
It may therefore be possible that, an improvement of \cref{thm:permutonLimit} (establishing the convergence of $\mu_n$ to $\mu$) 
with an additional concentration result could be used to prove Conjecture 5.2.6 in~\cite{Jacopo}.

\subsection{On guessing the candidate limit $\mu$}

Our proof of the convergence of $\mu_n$ to $\mu$ in \cref{ssec:proof_permuton_cv} 
provides no information on why~$\mu$ has been defined as it is. 
The purpose of this section is therefore to explain how we identified the candidate limit $\mu$. 
We emphasize that this is not necessary for the proof to come, and included only so that readers can understand why $\mu$ is the natural candidate for the limit, and how we found it. 

All along this section, we fix a value $\lambda>0$, and for each $n$, we let $\sigma_n$ denote a random \recbiased permutation of size $n$, for the parameter $\theta = \lambda n$. 

\paragraph{Equation of the curve}
The curve of equation $y=f_\lambda(x)$ is, informally speaking, describing the limit of the broken lines which interpolate between the records of $\sigma_n$. 
We start by presenting how the expression of $f_\lambda$ was derived. 

For $\sigma$ a permutation of size $n$ and $i\in[n]$, let $\lmax(\sigma,i)=\max\{\sigma(j):j\in[i]\}$ be the last (and largest) record occurring at or before position $i$. 
We are interested in
$\PP(\lmax(\sigma_n,i)=j)$ for $i\approx xn$ and $j\approx yn$, and $x,y\in(0,1)$. 
The idea is that the equation of the curve is obtained by finding the value $y$ that maximizes this probability for any given $x\in(0,1)$. 

To compute this probability, we decompose permutations of size $n$, 
and we examine the contribution of each part of this decomposition to their weight in the record-biased distribution. 

A permutation $\sigma$ of size $n$ such that $\lmax(\sigma,i)=j$ is decomposed as follows. 
\begin{itemize}
 \item First, observe that necessarily $j\geq i$. 
So, we first choose the images by $\sigma$ of $[i]$ in $[j]$ (there are $\binom{j-1}{i-1}$ possibilities, as $j$ belongs to this set). Once chosen, permuting them in all possible ways yields a total contribution to the weight of $\rfact\theta i$.  
 \item Then we identify the positions of the $j-i$ remaining images that are smaller than $j$, all occurring after position $i$. 
There are $\binom{n-i}{j-i}$ possible sets of preimages for them, and $(j-i)!$ ways to assign them; as none of these elements are records, the total weight remains unchanged. 
 \item The remaining images (all those of value greater than $j$) are assigned to the remaining positions, in any possible fashion. 
The total weight on all possible permutations of these remaining images is $\rfact\theta {n-j}$. 
\end{itemize}

Therefore, we have, for $j\geq i$,
\begin{equation}
\PP(\lmax(\sigma_n, i)=j) = \frac{{\rfact\theta i}{\rfact \theta{n-j}}}{\rfact\theta n}\binom{j-1}{i-1}
\binom{n-i}{j-i}(j-i)! = \frac{{\rfact\theta i}{\rfact \theta{n-j}}}{\rfact\theta n}\frac{i}{j}\binom{j}{i}
\binom{n-i}{j-i}(j-i)!. 
\label{eq:probaLmax_i=j}
\end{equation}

Now we proceed informally by using $i=xn$, $j=yn$ (recall also that $\theta = \lambda n$). 
We want to establish asymptotic estimates for all the factors occurring in the above expression of $\PP(\lmax(\sigma_n, i)=j)$. First, 
\[
{\rfact\theta i} =\frac{\Gamma(\lambda n+xn)}{\Gamma(\lambda n)}  =\frac{\lambda n}{(\lambda+x)n}\frac{\Gamma(\lambda n+xn+1)}{\Gamma(\lambda n+1)} 
\sim \sqrt{\frac{\lambda}{\lambda+x}} \left(\frac{(\lambda+x)^{\lambda+x}}{\lambda^{\lambda}}\right)^{n} n^{xn}
e^{-xn}
\]
Similarly, we find 
\begin{align*}
{\rfact\theta {n}} & \sim 
\sqrt{\frac{\lambda}{\lambda+1}} \left(\frac{(\lambda+1)^{\lambda+1}}{\lambda^{\lambda}}\right)^{n} n^{n}
e^{-n} \text{ and}\\
{\rfact\theta {n-j}} & \sim 
\sqrt{\frac{\lambda}{\lambda+1-y}} \left(\frac{(\lambda+1-y)^{\lambda+1-y}}{\lambda^{\lambda}}\right)^{n} n^{(1-y)n}
e^{-(1-y)n}.
\end{align*}
Hence
\[
 \frac{{\rfact\theta i}{\rfact \theta{n-j}}}{\rfact\theta n}
\sim \sqrt{\frac{(\lambda+x)(\lambda+1-y)}{\lambda(\lambda+1)}}
\left(\frac{(\lambda+x)^{\lambda+x}(\lambda+1-y)^{\lambda+1-y}}{\lambda^{\lambda}(\lambda+1)^{\lambda+1}}\right)^{n}
n^{(x-y)n}e^{-(x-y)n}.\]
For the factors which do not involve $\theta$, we have 
\[
\binom{n-i}{j-i}(j-i)! = \frac{\Gamma(n-i+1)}{\Gamma(n-j+1)}
\sim \sqrt{\frac{1-x}{1-y}} \left(\frac{(1-x)^{1-x}}{(1-y)^{1-y}}\right)^{n}e^{(x-y)n}n^{-(x-y)n}, 
\]
and 
\[
\frac{i}{j}\binom{j}{i} \sim \frac1{\sqrt{2\pi n}} \frac{x}{y} \sqrt{\frac{y}{x(y-x)}} \left(\frac{y^{y}}{x^{x}(y-x)^{y-x}}\right)^{n}.
\]
Putting everything together, we obtain 
\begin{equation}
\PP(\lmax(\sigma_n, i)=j) \sim \frac1{\sqrt{2\pi n}}  \sqrt{\frac{x(1-x)(\lambda+x)(\lambda+1-y)}{y(1-y)(y-x)\lambda(\lambda+1)}} F(x,y)^{n},
\label{eq:equivPLmax_informal}
\end{equation}
with
\begin{equation}
F(x,y) = \frac{y^{y}(1-x)^{1-x}(\lambda+x)^{\lambda+x}(\lambda+1-y)^{\lambda+1-y}}{x^{x}(y-x)^{y-x}(1-y)^{1-y}\lambda^{\lambda}(\lambda+1)^{\lambda+1}}.
\label{eq:defF}
\end{equation}
It is necessary that $F(x,y)\leq 1$, otherwise the probability would be greater than $1$ for $n$ sufficiently large. 
Even more, for any given $x$, we claim that there is at least one $y$ such that $F(x,y)=1$. 
Indeed, assuming otherwise, by continuity of $F$ we can make $\sum_{j}\PP(\lmax(\sigma_n, i)=j)$ as small as we want, which is not possible as it is equal to $1$. 

Therefore, the equation of the curve we are looking for is given by finding, for each $x$, the value of $y$ such that $F(x,y) = 1$, \emph{i.e.} which maximizes $F(x,y)$. 
Rather than solving $F(x,y) = 1$ or $\frac{d}{dy}F(x,y) = 0$, we can solve $\frac{d}{dy}\log F(x,y) = 0$, for which 
computations are easier, since  
\[
\frac{d}{dy}\log F(x,y) = \log\left(\frac{y(1-y)}{(y-x)(\lambda+1-y)}\right).
\]
Hence
\[
\frac{d}{dy}\log F(x,y) = 0 \Leftrightarrow y = \frac{x(\lambda+1)}{\lambda +x}, 
\]
which justifies our definition of the function $f_\lambda(x) = \frac{x(\lambda+1)}{\lambda +x}$ in the previous section. 

By contruction, it holds that 
\begin{equation}
\text{for all $x$, } F(x,f_\lambda(x))=1.
\label{eq:F(x,y)=1}
\end{equation}
It is also possible to check (\emph{e.g.} in Maple) from the explicit expression of $F(x,y)$ and of $f_\lambda(x)$ that the equality $F(x,f_\lambda(x))=1$ indeed holds for all $x$.

\paragraph{Density on the curve}
Next, we want to determine the density of the measure $\mu_c$ whose support is the curve of equation $y = f_\lambda(x)$. 
Informally, the total mass on this curve is the limit of the measure carried by all the records in the permutation $\sigma_n$. 
To compute it, let us first 
recall that, in a permutation of size $n$, $\frac{1}{n}$ is the measure associated with every point of $\sigma$ (to ensure that the total measure of all points is $1$). 
So, using \cref{lem:records}, we obtain that the total measure carried by all the records in $\sigma_n$ is 
\[
 \sum_{i=1}^n \frac{1}{n} \PP_{n}(\text{record at }i)  = \frac{1}{n} \sum_{i=1}^n \frac{\theta}{\theta  +  i - 1} =
\frac{1}{n} \sum_{i=1}^n \frac{\lambda}{\lambda  +  i/n - 1/n}.
\]
Up to the $1/n$ (that we neglect), this is a Riemann sum, and taking $n \to \infty$, 
the limit is 
$
\int_0^1 \frac{\lambda}{\lambda+x} dx.
$
If we were considering only records occurring between position $i \approx an$ and $j \approx bn$, 
the measure carried by these records would similarly be 
$
\int_a^b \frac{\lambda}{\lambda+x} dx.
$
This explains how the density $\frac{\lambda}{\lambda+x}$ of the measure $\mu_c$ was determined. 

\paragraph{Density below the curve}
Finally, we need to ``complete'' $\mu_c$ by a measure $\mu_u$ under the curve of equation $y=f_\lambda(x)$ 
in such a way that we obtain a permuton. In particular, we must ensure the uniform marginals property. 

By definition, 
\[
 \mu_c([a,b]\times[0,1]) = \int_a^b \frac{\lambda}{\lambda+x} dx. 
\]
Therefore, if we want to put a uniform measure of area density $a_\lambda$ under the curve, it needs to satisfy
\[
 \text{for all } a < b \in [0,1], \int_a^b \left( \frac{\lambda}{\lambda+x} + a_\lambda f_\lambda(x) \right) dx = b-a.  
\]
Recalling that $f_\lambda(x) = \frac{x(\lambda+1)}{\lambda +x}$, 
the condition above is satisfied when taking $a_{\lambda}=\frac{1}{\lambda +1}$. 

\subsection{Proof of the permuton limit}\label{ssec:proof_permuton_cv}

Throughout this section, we fix $\lambda >0$. 
Recall that, for each integer $n$, $\sigma_n$ denotes a random \recbiased permutation of size $n$, for the parameter $\theta = \lambda n$. 
Recall also that for $i\in[n]$, $\lmax(\sigma_n,i)=\max\{\sigma(j):j\in[i]\}$ denotes the last (and largest) record occurring at or before position $i$. 
Finally, recall that $\mu_n$ is the (random) permuton associated with $\sigma_n$, and that the permuton $\mu$ has been defined in \cref{sec:def_mu}. 

The rest of the paper is devoted to proving our main result in this section. 

\begin{theorem}
 \label{thm:permutonLimit}
%Let $\sigma_n$ be a random \recbiased permutation of size $n$, for the parameter $\theta = \lambda n$, and let $\mu_n$ be the corresponding (random) permuton. 
It holds that 
$\mu_n$ converges in probability to~$\mu$ (w.r.t. the weak topology).  
\end{theorem} 

\medskip

The first step of our proof is to show that $\lmax(\sigma_n,i)$ is concentrated around its typical value. This appears below as \cref{pro:large dev lmax}. 

\begin{proposition}
For every $\varepsilon>0$, there exists  $c(\varepsilon) \in(0,1)$ such that for $n$ sufficiently large we have:
\[ \text{for all } i \in\{1,\ldots,n\} \text{ such that } \varepsilon < f_\lambda(i/n) < 1-\varepsilon, \ 
\mathbb{P}\left(|\lmax(\sigma_n,i) - nf_\lambda(i/n)|>\varepsilon\,n\right) \leq c(\varepsilon)^n.
\]
\label{pro:large dev lmax}
\end{proposition}

Establishing this result essentially consists in transforming the estimate in \cref{eq:equivPLmax_informal} -- which is only informal (and qualitative) -- into a precise (and quantitative) statement. 
To this effect, we shall need tight upper and lower bounds on the rising factorials (\cref{lm:bounds gamma}). We then use them to establish upper and lower bounds on $\PP(\lmax(\sigma_n, i)=j)$ (\cref{lem:bornesPLmax(i)=j}). Both bounds are of the same order as the estimate in \cref{eq:equivPLmax_informal} up to polynomial factors. 

\begin{lemma}\label{lm:bounds gamma}
For any positive integer $n$ and any real $\alpha\geq 0$ such that $\alpha n$ is an integer, we have
\[
\frac1{c_\lambda}\left(\frac{(\lambda +\alpha)^{\lambda +\alpha }}{\lambda ^\lambda e^{\alpha}}\right)^n\,n^{\alpha n} \leq \rfact{(\lambda n)}{\alpha n} \leq c_\lambda \left(\frac{(\lambda +\alpha)^{\lambda +\alpha }}{\lambda ^\lambda e^{\alpha}}\right)^n\,n^{\alpha n},
\text{ with }c_\lambda = \sqrt{\frac{\lambda+1}\lambda}.
\]
\end{lemma}

\begin{proof}
We use the fact that for all $x>0$, we have~\cite[Thm 1.6]{Batir08}
\begin{equation}\label{eq:gamma}
x^x e^{-x} \sqrt{2\pi x} \leq \Gamma(x+1)\leq x^x e^{-x} \sqrt{2\pi(x+1)}.
\end{equation}

As $\rfact{\theta}i=\frac{\Gamma(\theta+i)}{\Gamma(\theta)}= \frac{\theta}{\theta +i} \frac{\Gamma(\theta+i+1)}{\Gamma(\theta+1)}$ we have, for $n\geq 1$
\begin{align*}
\rfact{(\lambda n)}{\alpha n}  & \leq  \frac{\lambda n}{\lambda n + \alpha n} \frac{((\lambda n +\alpha n)e^{-1})^{\lambda n+\alpha n} \sqrt{2\pi(\lambda n+\alpha n +1)}}{((\lambda n)e^{-1})^{\lambda n}\sqrt{2\pi \lambda n}}\\
&\leq  \frac{\lambda }{\lambda + \alpha } \left(\frac{(\lambda +\alpha)^{\lambda +\alpha }}{\lambda ^\lambda e^{\alpha}}\right)^n \frac{\sqrt{(\lambda +\alpha +1)}}{\sqrt{\lambda}}\,n^{\alpha n}\\
& \leq d_\lambda  \left(\frac{(\lambda +\alpha)^{\lambda +\alpha }}{\lambda ^\lambda e^{\alpha}}\right)^n\,n^{\alpha n},
\end{align*}
where $d_\lambda = \max_{\alpha\geq 0}\frac{\sqrt{\lambda(\lambda+\alpha+1)}}{\lambda+\alpha}
= \sqrt{\frac{\lambda+1}\lambda}.$
The lower bound works similarly:
\begin{align*}
\rfact{(\lambda n)}{\alpha n}  & \geq  \frac{\lambda n}{\lambda n + \alpha n} \frac{((\lambda n +\alpha n)e^{-1})^{\lambda n+\alpha n} \sqrt{2\pi(\lambda n+\alpha n)}}{((\lambda n)e^{-1})^{\lambda n}\sqrt{2\pi(\lambda n+1)}}\\
&\geq  \frac{\lambda }{\lambda + \alpha } \left(\frac{(\lambda +\alpha)^{\lambda +\alpha }}{\lambda ^\lambda e^{\alpha}}\right)^n \frac{\sqrt{\lambda +\alpha}}{\sqrt{\lambda+1}}\,n^{\alpha n}\\
& \geq e_\lambda  \left(\frac{(\lambda +\alpha)^{\lambda +\alpha }}{\lambda ^\lambda e^{\alpha}}\right)^n\,n^{\alpha n},
\end{align*}
where $e_\lambda = \min_{\alpha\geq0}\frac{\lambda}{\sqrt{\lambda+1}}\cdot\frac1{\sqrt{\lambda+\alpha}} = \sqrt{\frac{\lambda}{\lambda+1}}$. This concludes the proof.
\end{proof}

Recall from \cref{eq:defF} that $F(x,y) = \frac{y^{y}(1-x)^{1-x}(\lambda+x)^{\lambda+x}(\lambda+1-y)^{\lambda+1-y}}{x^{x}(y-x)^{y-x}(1-y)^{1-y}\lambda^{\lambda}(\lambda+1)^{\lambda+1}}$.

\begin{lemma}
For $n$ sufficiently large, and for all $1<i\leq j\leq n-1$, it holds that 
\begin{equation*}
n^{-5} \cdot F(\tfrac{i}n,\tfrac{j}n)^n\leq \mathbb{P}(\lmax(\sigma_n,i)=j) \leq n^3 \cdot F(\tfrac{i}n,\tfrac{j}n)^n.\end{equation*}
 \label{lem:bornesPLmax(i)=j}
\end{lemma}

\noindent Note that only the upper bound in the above proposition is needed later. However, obtaining the lower bound does not require much more work, so we record it here for completeness. 

\begin{proof}
Let $j$ be an integer such that  $i\leq j < n$. 
From \cref{eq:probaLmax_i=j}, we have
\[
\mathbb{P}(\lmax(\sigma_n,i)=j) = \frac{i}{j}\frac{\theta^{(i)}\theta^{(n-j)}}{\theta^{(n)}}\binom{j}{i}\binom{n-i}{j-i}(j-i)!.
\]
Let $\alpha = i/ n$ and $\beta =j/n$ (noting that this implies $\alpha \leq \beta$). From Lemma~\ref{lm:bounds gamma}, we have 
\begin{align*}
\frac{\theta^{(i)}\theta^{(n-j)}}{\theta^{(n)}} & \leq
c_\lambda^3 \left(\frac{(\lambda+\alpha)^{\lambda+\alpha}}{\lambda^\lambda e^\alpha}\right)^n
\left(\frac{(\lambda+1-\beta)^{\lambda+1-\beta}}{\lambda^\lambda e^{1-\beta}}\right)^n
\left(\frac{(\lambda+1)^{\lambda+1}}{\lambda^\lambda e}\right)^{-n} n^{(\alpha-\beta)n}\\
&\leq c^3_\lambda \left(\frac{(\lambda+\alpha)^{\lambda+\alpha}(\lambda+1-\beta)^{\lambda+1-\beta}}{\lambda^\lambda (\lambda+1)^{\lambda+1} e^{\alpha-\beta}}\right)^n n^{(\alpha-\beta)n}.
\end{align*}
And similarly,
\begin{align*}
\frac{\theta^{(i)}\theta^{(n-j)}}{\theta^{(n)}} & \geq \frac1{c^3_\lambda} \left(\frac{(\lambda+\alpha)^{\lambda+\alpha}(\lambda+1-\beta)^{\lambda+1-\beta}}{\lambda^\lambda (\lambda+1)^{\lambda+1} e^{\alpha-\beta}}\right)^n n^{(\alpha-\beta)n}.
\end{align*}
We now use the following (crude) Stirling bounds:
\begin{equation}\label{eq:stirling}
\forall n\geq 0,\ n^{n}e^{-n} \leq n!\quad \text{and}\quad \forall n\geq 1,\ n! \leq 3n^{n+1}e^{-n}.
\end{equation}
Hence
\[
\binom{j}{i} = \frac{j!}{i!(j-i)!}
\leq \frac{3j^{j+1}}{e^j}\frac{e^i}{i^i}\frac{e^{j-i}}{(j-i)^{j-i}}
\leq 3n \left(\frac{\beta^\beta}{\alpha^\alpha (\beta-\alpha)^{\beta-\alpha}}\right)^n,
\]
and
\[
\binom{n-i}{j-i} (j-i)! =  \frac{(n-i)!}{(n-j)!}
\leq \frac{3(n-i)^{n-i+1}}{e^{n-i}} \frac{e^{n-j}}{(n-j)^{n-j}}
\leq 3n \left(e^{\alpha-\beta} \frac{(1-\alpha)^{1-\alpha}}{(1-\beta)^{1-\beta}} \right)^n n^{(\beta-\alpha)n}.
\]
Putting all together, we get
\begin{align*}
\mathbb{P}(\lmax(\sigma_n,i)=j) & \leq  9n^2 c_\lambda^3 \frac{\alpha}{\beta}\left[\frac{(\lambda+\alpha)^{\lambda+\alpha}(\lambda+1-\beta)^{\lambda+1-\beta}}{\lambda^\lambda(\lambda+1)^{\lambda+1}}\right]^n
 \left(\frac{\beta^\beta}{\alpha^\alpha (\beta-\alpha)^{\beta-\alpha}}\right)^n\left( \frac{(1-\alpha)^{1-\alpha}}{(1-\beta)^{1-\beta}} \right)^n \\
& \leq 9n^2 c_\lambda^3 \frac{\alpha}{\beta} F(\alpha,\beta)^n \leq b_\lambda n^2 F(\alpha,\beta)^n,
\end{align*}
where $b_\lambda = 9 c_\lambda^3$ and $F$ is defined in \cref{eq:defF} as recalled above. 
 Similarly,
 \[
\binom{j}{i} = \frac{j!}{i!(j-i)!}\geq
\frac1{9n^2} \left(\frac{\beta^\beta}{\alpha^\alpha (\beta-\alpha)^{\beta-\alpha}}\right)^n
\text{ and }
\binom{n-i}{j-i} (j-i)! \geq \frac1{3n} \left(e^{\alpha-\beta} \frac{(1-\alpha)^{1-\alpha}}{(1-\beta)^{1-\beta}} \right)^n n^{(\beta-\alpha)n}.
\]
This yields
\[
\mathbb{P}(\lmax(\sigma_n,i)=j)  \geq
\frac1{27n^3c_\lambda^3}\frac{\alpha}{\beta} F(\alpha,\beta)^n \geq a_\lambda n^{-4} F(\alpha,\beta)^n \text{ for $n$ large enough},
\]
with $a_\lambda = \frac1{27c^3_\lambda}$. 
This completes the proof of \cref{lem:bornesPLmax(i)=j}.
\end{proof}

\begin{proof}[Proof of \cref{pro:large dev lmax}]
Given the upper bound of \cref{lem:bornesPLmax(i)=j}, we study $y\mapsto F(x,y)$ when $|x-y|>\varepsilon$, and we want an upper bound that is uniform in $x$. 
Rather than studying $F$ directly, we study the functions $H_x$ defined, for all $x \in (0,1)$, by $H_x(y) = \log F(x,y)$ for all $y \in (x,1)$ (which is computationally simpler). In other words, 
\begin{align*}
H_x(y) \ = \ & y \log(y) + (1-x) \log(1-x) + (\lambda +x) \log(\lambda +x) + (\lambda+1-y) \log(\lambda+1-y) \\
&- x \log(x) - (y-x)\log(y-x) -(1-y)\log(1-y) -\lambda \log(\lambda) -(\lambda+1)\log(\lambda +1). 
\end{align*}
Some properties of these functions, useful in the following, are immediately obtained. In particular, for any $x \in (0,1)$, it holds 
that  $H_x(f_\lambda(x)) =0$ (since $F(x,f_\lambda(x))=1$ as seen earlier), 
that $H_x$ is twice differentiable (w.r.t. $y$), 
that $\frac{d}{dy} H_x(f_\lambda(x)) =0$, and 
that $H_x$ is increasing on $(x,f_\lambda(x))$ and decreasing on $(f_\lambda(x),1)$. Also observe that, for $x<y<1$,
\begin{align*}
\frac{d^2}{dy^2} H_x(y) & = \frac1y - \frac1{1-y} - \frac1{y-x} + \frac1{\lambda+1-y}
= - \frac{x}{y(y-x)} - \frac{\lambda}{(1-y)(\lambda+1-y)}\\
& \leq - \frac{\lambda}{(1-y)(\lambda+1-y)} \leq \frac{-\lambda}{\lambda+1}.
\end{align*}

For the rest of the proof, fix $\varepsilon>0$. 
Assume $x \in (0,1)$ is such that $f_\lambda(x)\leq 1-\varepsilon$. 
Applying Taylor's theorem (with the Lagrange form of the remainder) yields that there exists $\zeta\in(x,1)$
such that
\[
 H_{x}\left(f_\lambda(x) +\varepsilon \right)
 =  H_{x}\left(f_\lambda(x)\right) + \varepsilon \frac{d}{dy}H_{x}\left(f_\lambda(x)\right) + \frac{\varepsilon^2}2 \frac{d^2}{dy^2}H_{x}\left(\zeta\right).
\]
The previous observations about $H_x$ therefore yield 
\[
 H_{x}\left(f_\lambda(x) +\varepsilon \right) = \frac{\varepsilon^2}2 \frac{d^2}{dy^2}H_{x}\left(\zeta\right) \leq \frac{-\lambda\varepsilon^2}{2(\lambda+1)}. 
\]
Similarly, assuming that $x \in (0,1)$ satisfies $f_\lambda(x)\geq \varepsilon$, we get 
\[
 H_{x}\left(f_\lambda(x) -\varepsilon \right) \leq \frac{-\lambda\varepsilon^2}{2(\lambda+1)}. 
\]
Using that $H_x$ is increasing until $y=f_\lambda(x)$ and then decreasing, 
it follows that, for $x \in (0,1)$ such that $\varepsilon \leq f_\lambda(x)\leq 1-\varepsilon$, for all $y \in (x,f_\lambda(x)-\varepsilon] \cup [f_\lambda(x)+\varepsilon,1)$, $H_{x}\left(y\right) \leq \frac{-\lambda\varepsilon^2}{2(\lambda+1)}$,  and thus 
\begin{equation}
F(x,y) \leq \exp\left(\frac{-\lambda\varepsilon^2}{2(\lambda+1)}\right),  \text{ a bound which does not depend on } x. 
\label{eq:borneSurF}
\end{equation}

We can now conclude the proof of \cref{pro:large dev lmax}. 
We have that, for $n$ sufficiently large, for all $i \in \{1, \dots, n\}$  such that  $\varepsilon < f_\lambda(i/n) < 1-\varepsilon$, 
\begin{align*}
\mathbb{P}(|\lmax(\sigma_n,i) & - nf_\lambda(i/n)|>\varepsilon\,n) \\
& = \sum_{j = 1}^{\lfloor nf_\lambda(i/n)-\varepsilon n\rfloor} \mathbb{P}(\lmax(\sigma_n,i) = j)
+\sum_{j = \lceil nf_\lambda(i/n)+\varepsilon n\rceil}^{n} \mathbb{P}(\lmax(\sigma_n,i) = j) \\
& \leq n \cdot n^3 \cdot \exp\left(\frac{-\lambda\varepsilon^2}{2(\lambda+1)}\right)^n, 
\end{align*}
where in the last line we used \cref{lem:bornesPLmax(i)=j} and \cref{eq:borneSurF}. This proves the announced statement that, for $n$ sufficiently large, $\mathbb{P}\left(|\lmax(\sigma_n,i) - nf_\lambda(i/n)|>\varepsilon\,n\right) \leq c(\varepsilon)^n$ for all $i \in \{1, \dots, n\}$ such that  $\varepsilon < f_\lambda(i/n) < 1-\varepsilon$, taking $c(\varepsilon)$ to be any value such that  $ \exp\left( \frac{-\lambda\varepsilon^2}{2(\lambda+1)}\right)<c(\varepsilon)<1$.
\end{proof}

\medskip

%%%%%%%%%%%%%%%%%%%%%%%%%%%%%%%

\begin{figure}[t]
\begin{center}
\begin{tikzpicture}[scale=3]
  \draw[] (0, 0) -- (1, 0) ;
  \draw[] (0, 0) -- (0, 1);
  \draw[dashed] (0,0) -- (1,1);
  \draw[dotted] (1,0) -- (1,1);
  \draw[dotted] (0,1) -- (1,1);
  \draw[red] (.3,.8) -- (.3,0) node[below]{$a$};
  \draw[red] (0.3,.8) -- (0,.8) node[left]{$b$};
  \draw[scale=1, domain=0:1, smooth, variable=\x, blue] plot ({\x}, {2*\x/(1+\x)});
  
  \node at (.5,-.3){\small tall rectangle};
  \end{tikzpicture}
  \qquad\qquad\qquad
\begin{tikzpicture}[scale=3]
  \draw[] (0, 0) -- (1, 0) ;
  \draw[] (0, 0) -- (0, 1);
  \draw[dashed] (0,0) -- (1,1);
  \draw[dotted] (1,0) -- (1,1);
  \draw[dotted] (0,1) -- (1,1);
  \draw[red] (.8,.3) -- (.8,0) node[below,anchor=north]{$a$};
  \draw[red] (0.8,.3) -- (0,.3) node[left]{$b$};
  \draw[scale=1, domain=0:1, smooth, variable=\x, blue] plot ({\x}, {2*\x/(1+\x)});
  
  \node at (.5,-.3){\small long rectangle};
  \end{tikzpicture}
  \end{center}
\caption{A rectangle starting in the bottom-left corner is either tall or long. \label{fig:long_and_tall_rectangles}}
\end{figure}
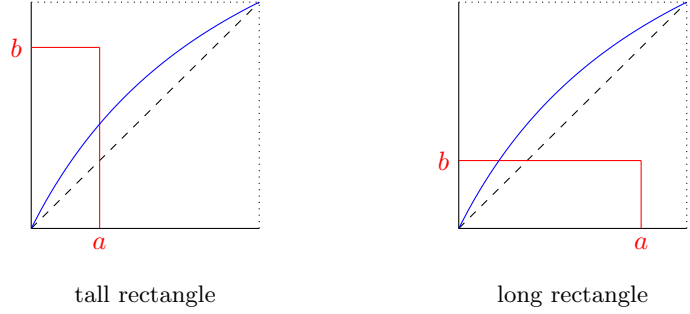

Now that the concentration of $\lmax(\sigma_n,i)$ has been established, we continue with our proof of \cref{thm:permutonLimit}. 
We intend to use \cref{lem:conv_permuton}, so that we will need to compare $\mu(R)$ and $\mu_n(R)$ for (grid-aligned) rectangles of the unit square having their bottom-left corner in $(0,0)$. 
This set of rectangles is naturally divided into two subsets, as illustrated by \cref{fig:long_and_tall_rectangles}: 
\begin{itemize}
\item the subset $\T$ of \emph{tall rectangles} consists of all rectangles $R=[0,a]\times[0,b]$ such that $b\geq f_\lambda(a)$; 
\item and the subset $\L$ of \emph{long rectangles} consists of all rectangles $R=[0,a]\times[0,b]$ such that $b<f_\lambda(a)$. 
\end{itemize}

%%%%%%%%%%%%%%%%%%%%%%%%%%%%%%%%%

With \cref{pro:large dev lmax}, it is easy to compare $\mu_n(R)$ and $\mu(R)$ for tall rectangles. 

\begin{proposition}
For every $\varepsilon>0$, there exists  $c(\varepsilon) \in(0,1)$ such that for $n$ sufficiently large we have, for any tall rectangle of the form $[0,i/n]\times[0,j/n]$ inside the unit square:
\[
\mathbb{P}\left(\left| \mu_n(R)-\mu(R)\right|>\varepsilon\right) \leq c(\varepsilon)^n.
\]
\label{prop:tall_rect}
\end{proposition}

\begin{proof}
Fix $\varepsilon >0$ and some $n$ large enough. 
Consider a tall rectangle $R=[0,i/n]\times[0,j/n]$ inside the unit square. 
Assume first that $\varepsilon < f_\lambda(i/n) < 1-\varepsilon$. 
From \cref{pro:large dev lmax}, it is enough to show that 
\[\mathbb{P}\left(\left| \mu_n(R)-\mu(R)\right|>\varepsilon\right) \leq
\mathbb{P}\left(|\lmax(\sigma_n,i) - nf_\lambda(i/n)|>\varepsilon\,n\right).\]
Because $R$ is a tall rectangle, by the uniform marginals property, we have $\mu(R) = i/n$. And because $\mu_n$ is a permuton, it holds that $\mu_n(R) \leq i/n$. Therefore, 
\[\mathbb{P}\left(\left| \mu_n(R)-\mu(R)\right|>\varepsilon\right) = \mathbb{P}\left(\mu_n(R)< i/n - \varepsilon\right).\]
Recall that $\mu_n$ is the permuton associated to the permutation $\sigma_n$ (of size $n$). 
Assuming $\mu_n(R)< i/n - \varepsilon$ means that 
$|\{k \leq i : \sigma_n(k) \leq j\}| < i - \varepsilon n$, that is to say that 
$|\{k \leq i : \sigma_n(k) > j\}| > \varepsilon n$. 
This last condition implies (by definition of $\lmax$) that 
$\lmax(\sigma_n,i) > j + n\varepsilon$. 
Since $R$ is a tall rectangle, we have $j/n \geq f_\lambda(i/n)$ so that $j + n\varepsilon \geq nf_\lambda(i/n) +\varepsilon n$. 
Summing things up, assuming that  $\mu_n(R)< i/n - \varepsilon$ implies that
$\lmax(\sigma_n,i) > nf_\lambda(i/n) +\varepsilon n$. 
This shows that 
\[
 \mathbb{P}\left(\mu_n(R)< i/n - \varepsilon\right) \leq 
\mathbb{P}\left(|\lmax(\sigma_n,i) - nf_\lambda(i/n)|>\varepsilon\,n\right)
\]
and therefore proves the proposition in the case where $\varepsilon < f_\lambda(i/n) < 1-\varepsilon$.

\medskip

If on the contrary $ f_\lambda(i/n) \leq \varepsilon$  or $f_\lambda(i/n) \geq 1-\varepsilon$, let us show that $|\mu_n(R) - \mu(R)| \leq \varepsilon$ deterministically, which will conclude the proof. 
First, as we already observed, $\mu_n(R) \leq \mu(R) = i/n$, so that we only need to show that $\mu(R) - \mu_n(R) \leq \varepsilon$. 

If $ f_\lambda(i/n) \leq \varepsilon$, then by definition of $f_\lambda$ we have $i/n \leq f_\lambda(i/n) \leq \varepsilon$ so that $\mu(R) - \mu_n(R) \leq \mu(R) = i/n \leq \varepsilon$. 

The final case is when $f_\lambda(i/n) \geq 1-\varepsilon$. Because $R$ is a tall rectangle, it holds that $j/n \geq f_\lambda(i/n)$, so that $\varepsilon \geq 1 - j/n$. Letting $S = [0,i/n]\times [0,1]$, we observe that $\mu_n(S) = i/n$ because of the uniform marginals property of permutons. Since $S$ is the disjoint union of $R$ and $T=[0,i/n]\times[j/n,1]$, it also holds that $\mu_n(R) = \mu_n(S) - \mu_n(T)$. Lastly, we note that $\mu_n(T) \leq \mu_n([0,1]\times[j/n,1]) = 1-j/n$ again by the uniform marginals property. These observations allow us to write that $\mu(R) - \mu_n(R) = i/n - \mu_n(S) + \mu_n(T) = \mu_n(T) \leq 1-j/n \leq \varepsilon$, concluding the proof. 
\end{proof}

\bigskip

We now turn to proving the analogue of \cref{prop:tall_rect} for long rectangles (see \cref{prop:long_rect} below). The proof is more involved, and we present it using some intermediate statements. 

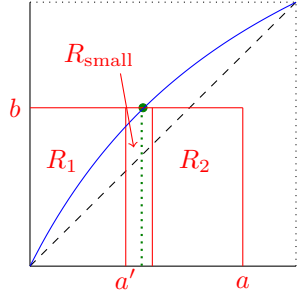
\begin{figure}[t]
\begin{center}
\begin{tikzpicture}[scale=3.5]
  \draw[] (0, 0) -- (1, 0) ;
  \draw[] (0, 0) -- (0, 1);
  \draw[dashed] (0,0) -- (1,1);
  \draw[dotted] (1,0) -- (1,1);
  \draw[dotted] (0,1) -- (1,1);

  \draw[red] (.36,.6) -- (.36,0) node[below,anchor=north,yshift=0.9mm]{$a'$};
  \draw[green!50!black,dotted,thick] (.42,.6) -- (.42,0);
  \draw[green!50!black,fill=green!50!black,] (.425,.6) circle(.015);
  \draw[red] (.8,.6) -- (.8,0) node[below,anchor=north,yshift=-0.2mm]{$a$};
  \draw[red] (0.8,.6) -- (0,.6) node[left]{$b$};
  \draw[scale=1, domain=0:1, smooth, variable=\x, blue] plot ({\x}, {2*\x/(1+\x)});
  \draw[red] (.46,.6) -- (.46,0);
  
  \node[red] at (.12,.4) {$R_1$};
  \node[red] at (.62,.4) {$R_2$};
  \node[red] at (.26,.8) {$R_{\text{small}}$};
  \draw[red,->] (.34,.73) -- (.39,.45);
\end{tikzpicture}
\end{center}
    \caption{Partition of a long rectangle $R$ into $R_1 \uplus R_{\text{small}} \uplus R_2$.}
    \label{fig:3rect}
\end{figure}

Recall that, for integers $i,j \leq n$, a rectangle $R = [0,i/n]\times[0,j/n]$ is long when $j/n < f_\lambda(i/n)$. 
For such a rectangle, we introduce the following notation:  $a=i/n$, $b=j/n$, $i'=\lfloor nf_\lambda^{-1}(b)\rfloor$, and $a'=i'/n$. 
We also define the rectangles $R_1 = [0,a']\times[0,b]$, $R_2 = [a'+1/n,a]\times [0,b]$ and $R_{\text{small}} = [a',a'+1/n]\times [0,b]$ (see \cref{fig:3rect}). Observe that $R_1$ is a tall rectangle, that $R_2$ entirely lies under the curve of $f_\lambda$, and that $R$ is the disjoint union of $R_1$, $R_2$ and $R_{\text{small}}$. In addition, because of the uniform marginals property of permutons, 
$\nu(R_{\text{small}}) \leq 1/n$ for any permuton $\nu$. 

In addition, for any $\varepsilon >0\in (0,1/2)$, and any long rectangle as above, we define, for all integers~$t$ such that $|t-nf_\lambda(a')|\leq \varepsilon^2 n$, the rectangle $R_t =[a'+1/n,a]\times [0,t/n]$. Observe that, because $t/n$ is close to $f_\lambda(a')$, $R_t$ does not differ much from $R_2$. 

Recall that $\sigma_n$ denotes a random \recbiased permutation of size $n$, for the parameter $\theta = \lambda n$, and that $\mu_n$ denotes the permuton associated with $\sigma_n$.
Recall also that for $i\in[n]$, $\lmax(\sigma_n,i)=\max\{\sigma(j):j\in[i]\}$ denotes the last (and largest) record occurring at or before position $i$. 

\medskip

Our strategy to prove \cref{prop:long_rect} is the following. First, in \cref{prop:union} we write $\mathbb{P}\left(|\mu_n(R)-\mu(R)|>30\,\varepsilon\right)$ as a sum of two terms, according to whether $\lmax(\sigma_n,i')$ is close or far from its typical value, and we establish an exponentially small bound for the latter case. Then, in 
\cref{prop:Rt}, we analyse the former case  to also prove that it is exponentially small, relying on a technical statement given in \cref{lem:aux}.

\begin{proposition}
\label{prop:union}
For any $\varepsilon\in (0,1/2)$,  there exists  $c(\varepsilon) \in(0,1)$ such that, for $n$ sufficiently large and for any long rectangle $R=[0,i/n]\times[0,j/n]$ inside the unit square, using the notation $i'=\lfloor nf_\lambda^{-1}(j/n)\rfloor$, and $a'=i'/n$ introduced above, we have: 
\begin{align*}
& \mathbb{P}\left(|\mu_n(R)-\mu(R)|>30\,\varepsilon\right)  \\ 
&  \qquad \leq \sum_{t \text{ s.t. } \atop {|t-nf_\lambda(a')|\leq \varepsilon^2 n}} \hspace{-0.6cm}\mathbb{P}\left(|\mu_n(R_t)-\mu(R_t)|> 8 \,\varepsilon \mid \lmax(\sigma_n,i')=t\right) \cdot \mathbb{P}(\lmax(\sigma_n,i')=t) + c(\varepsilon)^n.
\end{align*}
\end{proposition}

\begin{proof}
Let us first prove the statement when $f_\lambda(i'/n) <\varepsilon^2$. 
Recalling that $i'=\lfloor nf_\lambda^{-1}(j/n)\rfloor$ and that $f_\lambda$ is increasing, we 
write that 
\[
\tfrac{i'}{n} = \tfrac{1}{n}  \lfloor nf_\lambda^{-1}(j/n)\rfloor \geq \tfrac{1}{n} \left( nf_\lambda^{-1}(j/n) -1 \right) = f_\lambda^{-1}(j/n) -1/n, 
\]
so that $j/n \leq f_\lambda(\tfrac{i'+1}{n})$. 
By the mean value theorem, and since $\max_{x\in[0,1]} |f_\lambda'(x)| = f_\lambda'(0) = \tfrac{\lambda+1}\lambda$, under the assumption that $f_\lambda(i'/n) <\varepsilon^2$, we have, for $n$ large enough, 
\[
f_\lambda(\tfrac{i'+1}{n}) \leq f_\lambda(i'/n) + \tfrac{1}{n} \tfrac{\lambda+1}\lambda \leq 2 \varepsilon^2 \leq \varepsilon
\]
(where we used the assumption that $\varepsilon\in (0,1/2)$ for the last inequality). 
It follows that for any permuton $\nu$, $\nu(R) \leq j/n \leq \varepsilon$. Therefore, $|\mu_n(R)-\mu(R)| \leq 2 \varepsilon$, so that  $\mathbb{P}\left(|\mu_n(R)-\mu(R)|>30\,\varepsilon\right) =0$, which obviously implies our proposition in this particular case.

For the rest of the proof, let us then assume that $f_\lambda(i'/n) \geq \varepsilon^2$.
Writing $R$ as $R_1 \uplus R_{\text{small}} \uplus R_2$, we first derive that 
\begin{align*}
\mathbb{P}\left(|\mu_n(R)-\mu(R)|>30\,\varepsilon\right) & \leq 
\mathbb{P}\left(|\mu_n(R_1)-\mu(R_1)|>10\,\varepsilon\right)\\ & 
+ \mathbb{P}\left(|\mu_n(R_{\text{small}})-\mu(R_{\text{small}})|>10\,\varepsilon\right)\\ & 
+ \mathbb{P}\left(|\mu_n(R_2)-\mu(R_2)|>10\,\varepsilon\right).
\end{align*}
Because $R_1$ is a tall rectangle, we know from \cref{prop:tall_rect} that there exists $c_{\text{tall}}(\varepsilon) \in(0,1)$ independent of $R$ such that, for $n$ sufficiently large $\mathbb{P}\left(|\mu_n(R_1)-\mu(R_1)|>10\,\varepsilon\right) \leq c_{\text{tall}}(\varepsilon)^n$. 
In addition, $|\mu_n(R_{\text{small}})-\mu(R_{\text{small}})| \leq 2/n$; as a consequence, it holds that for $n$ sufficiently large, $\mathbb{P}\left(|\mu_n(R_{\text{small}})-\mu(R_{\text{small}})|>10\,\varepsilon\right) =0$. We therefore focus on $\mathbb{P}\left(|\mu_n(R_2)-\mu(R_2)|>10\,\varepsilon\right)$. 

Again, we distinguish a special case, this time when $f_\lambda(i'/n) > 1 - \varepsilon^2$. 
From $i'=\lfloor nf_\lambda^{-1}(j/n)\rfloor$ and the increasing property of $f_\lambda$, we first derive that $f_\lambda(i'/n) \leq j/n$, so that $1-j/n \leq \varepsilon^2$. 
We now define the rectangles $S=[\tfrac{i'+1}{n},\tfrac{i}{n}]\times[0,1]$ and $T=[\tfrac{i'+1}{n},\tfrac{i}{n}]\times[j/n,1]$, which satisfy that $S = R_2 \uplus T$. Note that, for any permuton $\nu$, by the uniform marginals property, we have $\nu(S) = \tfrac{i-i'-1}{n}$ and $\nu(T) \leq 1-j/n$. Therefore, 
\[
|\mu_n(R_2)-\mu(R_2)| = |\mu_n(S)-\mu_n(T) - \mu(S)+\mu(T)| = |\mu(T) - \mu_n(T)| \leq 2 \cdot (1-j/n) \leq 2\varepsilon^2\leq \varepsilon.
\] 
Consequently, $\mathbb{P}\left(|\mu_n(R_2)-\mu(R_2)|>10\,\varepsilon\right) = 0$ in this case.

We are now left with providing an upper bound on $\mathbb{P}\left(|\mu_n(R_2)-\mu(R_2)|>10\,\varepsilon\right)$, when the inequalities $\varepsilon^2 \leq f_\lambda(i'/n) \leq 1 - \varepsilon^2$ are satisfied.
We use the law of total probability to write that 
\begin{align*}
\mathbb{P}\left(|\mu_n(R_2)-\mu(R_2)|>10\,\varepsilon\right) & = \sum_{t=1}^n \mathbb{P}\left(|\mu_n(R_2)-\mu(R_2)|> 10 \,\varepsilon \mid \lmax(\sigma_n,i')=t\right) \cdot \mathbb{P}(\lmax(\sigma_n,i')=t) \\
& \leq \hspace{-0.6cm} \sum_{t \text{ s.t. } \atop {|t-nf_\lambda(a')|\leq \varepsilon^2 n}} \hspace{-0.6cm}\mathbb{P}\left(|\mu_n(R_2)-\mu(R_2)|> 10 \,\varepsilon \mid \lmax(\sigma_n,i')=t\right) \cdot \mathbb{P}(\lmax(\sigma_n,i')=t) \\
& +  \hspace{-0.6cm} \sum_{t \text{ s.t. } \atop {|t-nf_\lambda(a')| > \varepsilon^2 n}}  \hspace{-0.6cm} \mathbb{P}(\lmax(\sigma_n,i')=t).
\end{align*}

For the term $\sum_{t}\mathbb{P}(\lmax(\sigma_n,i')=t)$ with the index $t$ ranging on values such that $|t-nf_\lambda(a')| > \varepsilon^2 n$, we shall use \cref{pro:large dev lmax}. 
Indeed, since it holds that $\varepsilon^2 \leq f_\lambda(i'/n) \leq 1 - \varepsilon^2$, we can write that, for $n$ sufficiently large, 
\[
\hspace{-0.6cm} \sum_{t \text{ s.t. } \atop {|t-nf_\lambda(a')| > \varepsilon^2 n}}  \hspace{-0.6cm} \mathbb{P}(\lmax(\sigma_n,i')=t) = \mathbb{P}(|\lmax(\sigma_n,i')-nf_\lambda(i'/n)| > \varepsilon^2 n) \leq c(\varepsilon^2)^n
\]
for the constant $c(\varepsilon^2)$ given by \cref{pro:large dev lmax}. 
As a consequence, there exists $c_{2}(\varepsilon)$ such that, for $n$ large enough, 
\[
\hspace{-0.6cm} \sum_{t \text{ s.t. } \atop {|t-nf_\lambda(a')| > \varepsilon^2 n}}  \hspace{-0.6cm} \mathbb{P}(\lmax(\sigma_n,i')=t) \leq c_2(\varepsilon)^n.
\]

Putting things together, to conclude the proof, it is enough to show that, for any integer $t$ such that $|t-nf_\lambda(a')| \leq \varepsilon^2 n$, 
\[
\mathbb{P}\left(|\mu_n(R_2)-\mu(R_2)|> 10 \,\varepsilon \mid \lmax(\sigma_n,i')=t\right) \leq 
\mathbb{P}\left(|\mu_n(R_t)-\mu(R_t)|> 8 \,\varepsilon \mid \lmax(\sigma_n,i')=t\right). 
\]
Let us fix some integer $t$ such that $|t-nf_\lambda(a')|\leq \varepsilon^2 n$. 
Recall that $R_t = [a'+1/n,a]\times [0,t/n]$ and $R_2 = [a'+1/n,a]\times [0,b]$. 
Consequently, if $t/n \leq b$ then $R_2 = R_t \uplus [a'+1/n,a]\times[t/n,b]$ and otherwise $R_t = R_2 \uplus [a'+1/n,a]\times[b,t/n]$. Therefore, for any permuton $\nu$, we have
\[
|\nu(R_2)-\nu(R_t)| \leq  |b-t/n|. 
\]
Next, using the assumption that $|t-nf_\lambda(a')|\leq \varepsilon^2 n$ and the fact  that $f_\lambda(a')\leq b \leq f_\lambda(a'+1/n)$, we have
\[
|b-t/n| \leq |b-f_\lambda(a')| + \varepsilon^2 \leq |f_\lambda(a'+1/n)-f_\lambda(a')| + \varepsilon^2.
\]
Using again the mean value theorem, and recalling that $\max_{x\in[0,1]} |f_\lambda'(x)| = \frac{\lambda+1}\lambda$, we obtain that $|f_\lambda(a'+1/n)-f_\lambda(a')|\leq \frac{\lambda+1}{\lambda  n}$. Hence, recalling our assumption that $\varepsilon\in(0,\frac12)$,  for $n$ sufficiently large, and for any permuton $\nu$, 
\begin{equation}\label{eq:R2-Rt}
|\nu(R_2)-\nu(R_t)| \leq |b-t/n| \leq \frac{\lambda+1}{\lambda  n} + \varepsilon^2 \leq 2\,\varepsilon^2 \leq \varepsilon.
\end{equation}
Using the above for $\nu=\mu$ and $\nu = \mu_n$, it follows that 
\[
\left|\mu_n(R_2) - \mu(R_2)\right|>10\,\varepsilon \Rightarrow \left|\mu_n(R_t) - \mu(R_t)\right|> 8\,\varepsilon,
\]
and the above holds deterministically.
This proves that 
\[
\mathbb{P}\left(|\mu_n(R_2)-\mu(R_2)|> 10 \,\varepsilon \mid \lmax(\sigma_n,i')=t\right) \leq 
\mathbb{P}\left(|\mu_n(R_t)-\mu(R_t)|> 8 \,\varepsilon \mid \lmax(\sigma_n,i')=t\right),
\]
hence concluding the proof. 
\end{proof}

\begin{proposition}
\label{prop:Rt}
For any $\varepsilon\in (0,1/2)$, there exists  $c(\varepsilon) \in(0,1)$ such that, for $n$ sufficiently large and for any long rectangle $R=[0,i/n]\times[0,j/n]$ inside the unit square, using the notation $i'=\lfloor nf_\lambda^{-1}(j/n)\rfloor$, and $a'=i'/n$ introduced above, for any integer~$t$ such that $|t-nf_\lambda(a')|\leq \varepsilon^2 n$, the rectangle $R_t = [a'+1/n,a]\times [0,t/n]$ satisfies: 
\[
\mathbb{P}\left(|\mu_n(R_t)-\mu(R_t)|> 8 \,\varepsilon \mid \lmax(\sigma_n,i')=t\right) \leq c(\varepsilon)^n.
\]
\end{proposition}

To prove the above, we use a large deviation result, stated in the following lemma. Its proof follows the same structure as that of \cref{pro:large dev lmax} and is differed to the end of this section.

\begin{lemma}
For any positive integers $k$, $\ell$ and $n$, such that $\ell\leq n$ and $k\leq n$, 
let us denote by $X_{\ell,n,k}$ the random variable that counts the number of integers smaller than or equal to $\ell$ in a uniform random size-$k$ subset of $[n]$. 

For any $\varepsilon>0$, there exists $c(\varepsilon) \in(0,1)$ such that, for $n$ sufficiently large, for all $1\leq \ell, k \leq n$, 
\[
\mathbb{P}\left(\left|X_{\ell,n,k}-\frac{\ell k}{n} \right|>\varepsilon n\right) \leq c(\varepsilon)^n.
\] \label{lem:aux}
\end{lemma}

\begin{proof}[Proof of \cref{prop:Rt}]
We start by computing $\mu(R_t)$, or rather finding a good estimate for it. 
Recall from \cref{eq:R2-Rt} that, under the assumptions that $\varepsilon\in (0,1/2)$ and $|t-nf_\lambda(a')|\leq \varepsilon^2 n$, it holds that $|b-t/n|  \leq \varepsilon$. 
By definition, we also have that $f_\lambda(a'+1/n) \geq b$. 

In the case that $f_\lambda(a'+1/n) > t/n$, $\mu(R_t)$ is easily computed. Indeed, then $R_t$ lies entirely under the curve of $f_\lambda$, so that $\mu(R_t)$ is proportional to its area; namely, $\mu(R_t) = \tfrac{(a-(a'+1/n))t}{(\lambda+1)n}$. If on the contrary $f_\lambda(a'+1/n) \leq t/n$, then $R_t$ is not entirely under the curve of $f_\lambda$, but almost. More precisely, in this case we have  $0\leq t/n - f_\lambda(a'+1/n) \leq t/n -b \leq \varepsilon$, from which we derive (using again the uniform marginal property of permutons) that $\mu(R_t) \leq \tfrac{(a-(a'+1/n))t}{(\lambda+1)n} + \varepsilon$. Therefore, in both cases, it holds that 
\[
\left|\mu(R_t) - \tfrac{(a-(a'+1/n))t}{(\lambda+1)n} \right|
\leq \varepsilon.
\]
As a consequence, \cref{prop:Rt} will follow if we show that 
there exists $c(\varepsilon)\in(0,1)$ that only depends on $\varepsilon$ such that for $n$ sufficiently large, it holds that 
\begin{equation}\label{eq:last reduction}
\mathbb{P}\left(\left|\mu_n(R_t) - \frac{(i-i'-1)t}{(\lambda+1)n^2} \right| > 7\,\varepsilon\mid\lmax(\sigma_n,i')=t\right) \leq c(\varepsilon)^n.
\end{equation}
To prove \cref{eq:last reduction}, we express $\mu_n(R_t)$ by a quantity defined from the permutation $\sigma_n$ directly. Indeed, $R_t$ has its sides supported by the $n \times n$ grid inside the unit square. This implies that $\mu_n(R_t) = \tfrac{1}{n} \cdot \rho_{i',i}^t(\sigma_n)$ where, for any permutation $\sigma$ of size $n$, we define $\rho_{i',i}^t(\sigma)$ as the number of elements $x$ of $\{i'+1, i'+2,\ldots, i-1\}$ such that $\sigma(x) \leq t$. \cref{eq:last reduction} can therefore be rewritten as 
\begin{equation}\label{eq:perm1}
\mathbb{P}\left(\left|\rho_{i',i}^t(\sigma_n)-\frac{t(i-i'-1)}{(\lambda+1)n}\right|>7\,\varepsilon\,n\mid  \lmax(\sigma_n,i')=t\right) \leq c(\varepsilon)^n.
\end{equation}

We now describe the structure of permutations of size $n$ satisfying the condition $\lmax(\sigma,i')=t$, which is equivalently written as $\max_{x\in[i']}\sigma(x) = t$. Such a permutation is entirely characterized by the data of 
\begin{itemize}
    \item a subset $T \in [t]$ of size $i'$ containing $t$, together with a bijective mapping $\phi_{\text{left}}$ from $[i']$ to~$T$ (this describes $\sigma([i']) =T$ and the permutation of these values occurring in the prefix of length $i'$ of $\sigma$); 
    \item a subset $A$ of size $t-i'$ of $\{i'+1, \ldots, n\}$, together with a bijective mapping $\phi_{\text{right,small}}$ from~$A$ to $S = [t] \setminus T$ (this describes the points of $\sigma$ occurring after position $i'$ that have a value less than $t$, as well as their images by $\sigma$); 
    \item a bijective mapping $\phi_{\text{right,large}}$ from $B = \{i'+1, \ldots, n\} \setminus A$ to $\{t+1, \ldots, n\}$ (this describes~$\sigma$ restricted to positions larger than $i'$ and values larger than $t$). 
\end{itemize}
Observe that, for a permutation $\sigma$ as above, $\rho_{i',i}^t(\sigma) = |A \cap [i'+1, i-1]|$. 
In addition, the number of records of a permutation $\sigma$ as above depends only on $\phi_{\text{left}}$ and $\phi_{\text{right,large}}$. In particular, for~$\sigma_n$ a \recbiased permutation of size $n$, and under the condition that  $\lmax(\sigma_n,i')=t$, the set~$A$ defined above is a uniform random subset of size $t-i'$ of $\{i'+1, \ldots, n\}$. 
As a consequence, using the notation of \cref{lem:aux}, $\rho_{i',i}^t(\sigma_n)$ is equal in distribution to $X_{i-1-i', n-i', t-i'}$ under the conditioning $\lmax(\sigma_n,i')=t$. 
The strategy to prove \cref{eq:perm1} is then to use \cref{lem:aux}, but it only applies when $n-i'$ tends to infinity. We therefore start by discarding the cases where $i'$ is too large. 

First observe that, for $\sigma$ a permutation of size $n$ such that $\lmax(\sigma,i')=t$, by definition we have $\rho_{i',i}^t(\sigma) \leq n-i'$, from which we deduce 
\[
\left|\rho_{i',i}^t(\sigma) - \frac{t(i-i'-1)}{(\lambda+1)n}\right|
\leq \rho_{i',i}^t(\sigma) +  \frac{n(n-i')}{(\lambda+1)n}
\leq \frac{\lambda+2}{\lambda+1}(n-i') \leq 2(n-i').
\]
Therefore, if $i'\geq n-3\varepsilon n$, then
$\mathbb{P}\left(\left|\rho_{i',i}^t(\sigma_n)-\frac{t(i-i'-1)}{(\lambda+1)n}\right|>7\,\varepsilon\,n\mid  \lmax(\sigma_n,i')=t\right) = 0$, proving our claim in this case. For the rest of the proof, let us then assume that $i'< n-3\varepsilon n$, that is to say $n-i' > 3\varepsilon n$. 
In particular, for any fixed $\varepsilon \in (0,1/2)$,  $n-i'$ tends to infinity when $n$ does. 
Applying \cref{lem:aux} (with $4\varepsilon$ instead of $\varepsilon$), we have that there exists $c_X(\varepsilon) \in(0,1)$ depending only on $\varepsilon$ such that, for $n$ sufficiently large, 
\begin{equation*}
\mathbb{P}\left(\left|\rho_{i',i}^t(\sigma_n)-\frac{(t-i')(i-i'-1)}{n-i'}\right|>4\,\varepsilon\,(n-i')\mid  \lmax(\sigma_n,i')=t\right) \leq c_X(\varepsilon)^{n-i'}\leq c_X(\varepsilon)^{3\varepsilon n}, 
\end{equation*}
which implies 
\begin{equation}
\mathbb{P}\left(\left|\rho_{i',i}^t(\sigma_n)-\frac{(t-i')(i-i'-1)}{n-i'}\right|>4\,\varepsilon\,n\mid  \lmax(\sigma_n,i')=t\right)\leq c_X(\varepsilon)^{3\varepsilon n}. 
\label{eq:presque}
\end{equation}
This resembles \cref{eq:perm1} which we want to prove: 
the differences are that the fraction $\frac{(t-i')(i-i'-1)}{n-i'}$ in \cref{eq:presque} is replaced by $\frac{t(i-i'-1)}{(\lambda+1)n}$ in \cref{eq:perm1}, and that we allow ourselves $7 \, \varepsilon \, n$ in \cref{eq:perm1} instead of $4 \, \varepsilon \, n$ in \cref{eq:presque}. 

To prove \cref{eq:perm1} from \cref{eq:presque}, we shall use small upper bounds of several differences, which we now give. Recall from \cref{eq:R2-Rt} that $|t/n - b| \leq 2 \varepsilon^2$ for $n$ large enough. Therefore (recalling also that $\varepsilon \in (0,1/2$)), for $n$ sufficiently large, it holds that 
\begin{equation}
\left|\frac{t(i-i'-1)}{(\lambda+1)n}-\frac{b(i-i'-1)}{(\lambda+1)}\right|
= \left|\frac{t}n - b\right|\cdot \frac{i-i'-1}{\lambda+1} \leq 2\,\varepsilon^2 n \leq \varepsilon n.
\label{eq:presque_1}
\end{equation}
In addition, using that $n-i'\geq 3\,\varepsilon n$, we also have, for $n$ sufficiently large
\begin{equation}
\left|\frac{(t-i')(i-i'-1)}{n-i'}-\frac{(bn-i')(i-i'-1)}{n-i'}\right| =n\cdot  \left|\frac{t}n-b\right|\cdot\frac{i-i'-1}{n-i'}\leq n\cdot 2\,\varepsilon^2\frac{n}{3\,\varepsilon\,n}\leq \varepsilon n.
\label{eq:presque_2}
\end{equation}
For the other two estimates that we need, it is convenient to introduce the notation $\beta = f_{\lambda}^{-1}(b)$. Because $i'$ is defined as $i'=\lfloor nf_\lambda^{-1}(b)\rfloor$, we have $i'/n \leq \beta \leq (i'+1)/n$. In addition, since $i'< n-3\varepsilon n$, it holds that $\beta \leq 1- 3 \varepsilon + 1/n \leq 1- 2 \varepsilon$ for $n$ large enough. Consequently, for $n$ sufficiently large, we have 
\[
\frac{bn-i'}{n-i'} = \frac{b-i'/n}{1-i'/n} \leq \frac{b-\beta +1/n}{1-\beta} = \frac{b-\beta}{1-\beta} +\frac{1}{1-\beta} \frac{1}{n} \leq \frac{b-\beta}{1-\beta} +\frac{1}{2 \varepsilon} \frac{1}{n}\leq \frac{b-\beta}{1-\beta} +\varepsilon.
\]
Therefore, for $n$ sufficiently large, 
\begin{equation}
\left|\frac{(bn-i')(i-i'-1)}{n-i'} - \frac{(b-\beta)(i-i'-1)}{1-\beta} \right| \leq \varepsilon \, n.
\label{eq:presque_4}
\end{equation}
Finally, since $b = f_\lambda(\beta) = \frac{\beta(\lambda+1)}{\lambda+\beta}$, we have 
\begin{equation}
\frac{b}{\lambda+1} = \frac{b- \beta}{1-\beta}. 
\label{eq:presque_3}
\end{equation}
Indeed, this can for example be derived from 
\begin{align*}
b = \frac{\beta(\lambda+1)}{\lambda+\beta}& \Leftrightarrow \frac{b}{\lambda+1}=\frac{\beta}{\lambda+\beta} \ \Leftrightarrow \ b(\lambda+\beta) = \beta(\lambda+1)\ \Leftrightarrow \ (b-\beta)\lambda= \beta-b\beta\\
&  \Leftrightarrow (b-\beta)(\lambda+1)= \beta-b\beta + (b-\beta) = (1-\beta)b \ \Leftrightarrow \ \frac{b}{\lambda+1} = \frac{b-\beta}{1-\beta}.
\end{align*}

Combining \cref{eq:presque_1} to \cref{eq:presque_3}, we derive that, for $n$ sufficiently large, 
\begin{align*}
\left|\rho_{i',i}^t(\sigma_n)-\frac{t(i-i'-1)}{(\lambda+1)n}\right|>7\,\varepsilon\,n
& \Rightarrow \left|\rho_{i',i}^t(\sigma_n)-\frac{b(i-i'-1)}{\lambda+1}\right|>6\,\varepsilon\,n \qquad \text{from \cref{eq:presque_1}} \\
& \Rightarrow \left|\rho_{i',i}^t(\sigma_n)-\frac{(b-\beta)(i-i'-1)}{1-\beta}\right|>6\,\varepsilon\,n \qquad \text{from \cref{eq:presque_3}}\\
& \Rightarrow \left|\rho_{i',i}^t(\sigma_n)-\frac{(bn-i')(i-i'-1)}{n-i'}\right|>5\,\varepsilon\,n \qquad \text{from \cref{eq:presque_4}}\\
& \Rightarrow \left|\rho_{i',i}^t(\sigma_n)-\frac{(t-i')(i-i'-1)}{n-i'}\right|>4\,\varepsilon\,n \qquad \text{from \cref{eq:presque_2}}. 
\end{align*}
As a consequence, for $n$ sufficiently large, it holds that 
\begin{align*}
& \mathbb{P}\left(\left|\rho_{i',i}^t(\sigma_n)-\frac{t(i-i'-1)}{(\lambda+1)n}\right|>7\,\varepsilon\,n\mid  \lmax(\sigma_n,i')=t\right) \\
& \qquad  \leq \mathbb{P}\left(\left|\rho_{i',i}^t(\sigma_n)-\frac{(t-i')(i-i'-1)}{n-i'}\right|>4\,\varepsilon\,n\mid  \lmax(\sigma_n,i')=t\right),
\end{align*}
and our proposition follows from the above and \cref{eq:presque}.
\end{proof}

With \cref{prop:union,prop:Rt} in our hands, the desired large deviation result for long rectangles is easily derived. 

\begin{proposition}
For every $\varepsilon \in (0,1/2)$, there exists  $c(\varepsilon) \in(0,1)$ such that for $n$ sufficiently large we have, for any long rectangle $R=[0,i/n]\times[0,j/n]$ inside the unit square:
\[
\mathbb{P}\left(\left| \mu_n(R)-\mu(R)\right|>30\varepsilon\right) \leq c(\varepsilon)^n.
\]
\label{prop:long_rect}
\end{proposition}

\begin{proof}
From \cref{prop:union}, we know that there exists $c_{\ref{prop:union}}(\varepsilon)$ such that for any long rectangle $R=[0,i/n]\times[0,j/n]$ inside the unit square, for $n$ large enough, 
\begin{align*}
& \mathbb{P}\left(|\mu_n(R)-\mu(R)|>30\,\varepsilon\right)  \\ 
&  \qquad \leq \sum_{t \text{ s.t. } \atop {|t-nf_\lambda(a')|\leq \varepsilon^2 n}} \hspace{-0.6cm}\mathbb{P}\left(|\mu_n(R_t)-\mu(R_t)|> 8 \,\varepsilon \mid \lmax(\sigma_n,i')=t\right) \cdot \mathbb{P}(\lmax(\sigma_n,i')=t) + c_{\ref{prop:union}}(\varepsilon)^n.
\end{align*}
Then, \cref{prop:Rt} gives us a constant $c_{\ref{prop:Rt}}(\varepsilon)$ such that, for any long rectangle $R=[0,i/n]\times[0,j/n]$ inside the unit square, and for any $t$ such that $|t-nf_\lambda(a')|\leq \varepsilon^2 n$, the rectangle $R_t$ satisfies that, for $n$ large enough,  
\[
\mathbb{P}\left(|\mu_n(R_t)-\mu(R_t)|> 8 \,\varepsilon \mid \lmax(\sigma_n,i')=t\right) \leq c_{\ref{prop:Rt}}(\varepsilon)^n.
\]
It follows that there exists $c(\varepsilon)$ such that for $n$ large enough,
\[
\mathbb{P}\left(|\mu_n(R)-\mu(R)|>30\,\varepsilon\right) \leq n c_{\ref{prop:Rt}}(\varepsilon)^n + c_{\ref{prop:union}}(\varepsilon)^n \leq c(\varepsilon)^n.\qedhere 
\]
\end{proof}

Combining \cref{prop:tall_rect,prop:long_rect} (large deviation results for tall and long rectangles), we can conclude the proof of the convergence of $\mu_n$ to $\mu$, as stated in \cref{thm:permutonLimit}.

\begin{proof}[Proof of \cref{thm:permutonLimit}]
From \cref{lem:conv_permuton}, it is enough to show that, for every $\varepsilon \in (0,1/2)$, we have $\PP(d_{\cdot /n} (\mu_n, \mu)>\varepsilon) \to 0$ as $n \to \infty$. Let us recall that the distance $d_{\cdot /n}$ is defined by $d_{\cdot /n}(\mu , \nu) = \sup (|\mu(R)-\nu(R)|)$
where $R$ ranges over the family $\mathcal{R}_n$ of rectangles of the form $[0,i/n]\times[0,j/n]$, for $0\leq i,j \leq n$. 
Using the union bound, we obtain 
\[
\PP(d_{\cdot /n}(\mu_n , \mu) >\varepsilon) \leq \sum_{R \in \mathcal{R}_n} \PP(|\mu_n(R)-\mu(R)| > \varepsilon). 
\]
Since there are only $\mathcal{O}(n^2)$ rectangles in $\mathcal{R}_n$, all of which are either long or tall, \cref{prop:tall_rect,prop:long_rect} ensure that $\PP(d_{\cdot /n}(\mu_n , \mu) >\varepsilon)$ indeed tends to $0$ as $n$ tends to infinity (and even at exponential speed). 
\end{proof}

We finish by giving the proof of \cref{lem:aux}, as announced above.
\begin{proof}[Proof of \cref{lem:aux}]
We first shrink the range of useful values for our parameters, as we are only interested in bounding from above $\mathbb{P}(|X_{\ell,n,k}-\frac{\ell k}n|>\varepsilon n)$. Since $X_{\ell,n,k}\leq \ell\leq n$ and $X_{\ell,n,k}\leq k\leq n$, if either $k$ or $\ell$ is smaller than $\varepsilon n$, then 
$X_{\ell,n,k}-\frac{\ell k}n \leq X_{\ell,n,k}\leq \varepsilon n$ and 
$\frac{\ell k}n - X_{\ell,n,k} \leq \frac{\ell k}n \leq \varepsilon n$. Hence, in these cases, 
$\mathbb{P}(|X_{\ell,n,k}-\frac{\ell k}n|>\varepsilon n) = 0$. 
 
 Assume now that $\ell \geq n-\varepsilon n$. Then, out of $k$ random elements, at most $\varepsilon n$ are not in $[\ell]$. Therefore $X_{\ell,n,k}\geq k-\varepsilon n$, yielding
 $\frac{k\ell}n-X_{\ell,n,k}\leq k -k +\varepsilon n\leq \varepsilon n$. Also, as $X_{\ell,n,k}\leq k$ we have $X_{\ell,n,k}-\frac{k\ell}n\leq k - \frac{k\ell}n \leq k - \frac{k(n-\varepsilon n)}n \leq \varepsilon k\leq \varepsilon n$. So in this case also, $\mathbb{P}(|X_{\ell,n,k}-\frac{\ell k}n|>\varepsilon n) = 0$. We can similarly prove that $k \geq n-\varepsilon n$ implies $\mathbb{P}(|X_{\ell,n,k}-\frac{\ell k}n|>\varepsilon n) = 0$. 
 Thus, we only consider the cases where 
$\varepsilon n \leq k\leq n-\varepsilon n$, and  $\varepsilon n \leq \ell\leq n-\varepsilon n$  in the sequel. 

\medskip
 
Clearly, $\mathbb{P}(X_{\ell,n,k} = r) =  \frac{\binom{\ell}{r}\binom{n-\ell}{k-r}}{\binom{n}{k}}$. This probability is $0$ if $r > \ell$ or $ r>k$ or $k-r>n-\ell$, so we can assume that $\ell \geq r$, $k \geq r$ and $r \geq k+\ell -n$ from now on. Using the (crude) Stirling bounds of \cref{eq:stirling} we have for any positive integers $a<b$:
\[
 \frac{a^a}{9a^2\cdot b^b(a-b)^{a-b}} \leq \binom{a}{b} \leq \frac{3a\cdot a^a}{b^b(a-b)^{a-b}},
\]
so that
\[
\frac{\binom{\ell}{r}\binom{n-\ell}{k-r}}{\binom{n}{k}}
\leq 81n^4 \frac{\ell^\ell(n-\ell)^{n-\ell}k^k(n-k)^{n-k}}{r^r(\ell-r)^{\ell-r}(k-r)^{k-r}(n-\ell-k+r)^{n-\ell-k+r}n^n}. 
\]
We define the function $\kappa$ on triples $(x,y,z) \in [0,1]^3$ such that $x-z\geq0$, $y-z\geq 0$ and $1-x-y+z \geq 0$ by \[
\kappa(x,y,z) = \frac{x^x(1-x)^{1-x}y^y(1-y)^{1-y}}{z^z(x-z)^{x-z}(y-z)^{y-z}(1-x-y+z)^{1-x-y+z}},
\]
with the convention that $0^0 =1$. 
Setting $x=\ell/n$, $y=k/n$ and $z=r/n$, the previous inequality can be rewritten as 
\[
\frac{\binom{\ell}{r}\binom{n-\ell}{k-r}}{\binom{n}{k}}
\leq 81\cdot n^4  \cdot\kappa\left(\frac{\ell}n,\frac{k}n,\frac{r}n\right)^n.
\]
A lower bound can be derived similarly. Putting things together yields that, for $n$ sufficiently large,
\begin{equation}\label{eq:kappa}
\frac1{n^6}\cdot \kappa\left(\frac{\ell}n,\frac{k}n,\frac{r}n\right)^n
\leq
\mathbb{P}(X_{\ell,n,k} = r) \leq n^5\cdot\kappa\left(\frac{\ell}n,\frac{k}n,\frac{r}n\right)^n.
\end{equation}
Using the assumptions 
$\varepsilon n \leq k, \ell \leq n - \varepsilon n$
%and $r \geq k+\ell -n$ 
described earlier in the proof, we are only interested in $\kappa(x,y,z)$ for 
$x,y\in [\varepsilon,1-\varepsilon]$. % and $z \geq \max\{0,x+y-1\}$. 
%\task{En fait, la définition de $\kappa$ suppose que $z \geq \max\{0,x+y-1\}$. Enlever d'ici cette condition ? }

For given $x,y\in [\varepsilon, 1-\varepsilon]$ we consider the real-valued map $\kappa_{x,y}: z\mapsto \kappa(x,y,z)$
defined on the (nonempty) interval $ I=[\max\{0,x+y-1\},\min\{x,y\}]$. This map is continuous, and infinitely differentiable on $\mathring{I}=(\max\{0,x+y-1\},\min\{x,y\})$. Simple computations yield:
\[
\frac{d}{dz}\log \kappa_{x,y}(z) = \log\frac{(x-z)(y-z)}{z(1-x-y+z)}.\]
Observe that
\[
\frac{d}{dz}\log \kappa_{x,y}(z) = 0 \Leftrightarrow (x-z)(y-z) = z(1-x-y+z)
\Leftrightarrow z=xy.
\]
Observe also that $xy\in \mathring{I}$. Moreover $\frac{d}{dz}\log \kappa_{x,y}(z)$ is positive for $z< xy$ and negative for $z>xy$. Therefore, $\log \kappa_{x,y}(z)$, and thus $\kappa_{x,y}(z)$, admits a unique maximum
at $z=xy$. Let $K=\kappa_{x,y}(xy)$. Because of the lower bound of \cref{eq:kappa}, 
we cannot have $K>1$, otherwise, by continuity, $\kappa_{x,y}$ would be greater than $1$ in an interval around $xy$, and we could build probabilities greater than $1$ for well chosen-values of the parameters. Hence $K\leq 1$.

Observe that on $\mathring{I}$, 
\[
\frac{d^2}{dz^2}\log \kappa_{x,y}(z) = -\frac1z - \frac1{x-z} - \frac1{y-z} - \frac1{1-x-y+z}
\leq - \frac1{x-z} \leq -\frac1{1-\varepsilon}.
\]
If $xy-\varepsilon\in\mathring{I}$ then, by Taylor's theorem (with the Lagrange form of the remainder), there exists $\zeta\in\mathring{I}$ such that
\[
\log \kappa_{x,y}(xy-\varepsilon) = 
\log \kappa_{x,y}(xy) + \frac{\varepsilon^2}{2} \frac{d^2}{dz^2}\log \kappa_{x,y}(\zeta)
\leq \log \kappa_{x,y}(xy) -\frac{\varepsilon^2}{2(1-\varepsilon)}.
\]
Therefore, as $\kappa_{x,y}(xy)=K\leq 1$ we have

\[
 \kappa_{x,y}(xy-\varepsilon) \leq K \exp\left( -\frac{\varepsilon^2}{2(1-\varepsilon)}\right) \leq \exp\left( -\frac{\varepsilon^2}{2(1-\varepsilon)}\right).
 \]
 Similarly, we also have, if $xy+\varepsilon\in\mathring{I}$,
\[
 \kappa_{x,y}(xy+\varepsilon) \leq \exp\left( -\frac{\varepsilon^2}{2(1-\varepsilon)}\right).
 \]
Because $\kappa_{x,y}$ is continuous, increasing until $z=xy$, then decreasing, this yields that:
\begin{equation}\label{eq:deviation kappa}
\forall z\in I\setminus (xy-\varepsilon,xy+\varepsilon),\ \kappa_{x,y}(z) \leq  \exp\left( -\frac{\varepsilon^2}{2(1-\varepsilon)}\right).
\end{equation}
We can now conclude the proof as follows. 
\begin{align*}
\mathbb{P}\left(\left| X_{\ell,n,k} - \frac{\ell k}{n}\right| > \varepsilon n\right) 
& =
\sum_{r<  \frac{\ell k}{n} - \varepsilon n} \mathbb{P}(X_{\ell,n,k} =r) + \sum_{r>  \frac{\ell k}{n} + \varepsilon n} \mathbb{P}(X_{\ell,n,k} =r) \\
& \leq n^5 \sum_{r<  \frac{\ell k}{n} - \varepsilon n} \kappa\left(\frac{\ell}n,\frac{k}n,\frac{r}n\right)^n
+  n^5 \sum_{r> \frac{\ell k}{n} + \varepsilon n} \kappa\left(\frac{\ell}n,\frac{k}n,\frac{r}n\right)^n\\
& \leq n^5 \sum_{r<  \frac{\ell k}{n} - \varepsilon n}  \exp\left( -\frac{n\varepsilon^2}{2(1-\varepsilon)}\right) + n^5 \sum_{r> \frac{\ell k}{n}  +\varepsilon n}  \exp\left( -\frac{n\varepsilon^2}{2(1-\varepsilon)}\right) \\
& \leq n^6 \exp\left( -\frac{n\varepsilon^2}{2(1-\varepsilon)}\right).
\end{align*}
So choosing any value $c(\varepsilon)$ such that $\exp\left( -\frac{\varepsilon^2}{2(1-\varepsilon)}\right) < c(\varepsilon) < 1$ proves the announced statement.
\end{proof}

\paragraph{Acknowledgments} 
We are grateful to Jacopo Borga and Valentin Féray for helpful discussions, comments and pointers at various stages of the preparation of this paper. 
We thank Rados\l aw Adamczak, Cyril Banderier, Sergi Elizalde and Svante Janson for their suggestions following presentations of the present work at AofA 2021 and at the Dagstuhl workshop on \emph{Pattern Avoidance, Statistical Mechanics and Computational Complexity} in March 2023. 

\bibliographystyle{plain}
\bibliography{versionlongue}

\end{document}